\documentclass
[
    a4paper,
    DIV=11,
    abstracton
]
{scrartcl}

\usepackage
{
    graphicx,
    amssymb,
    amsmath,
    amsthm,
    dsfont, % for 1 vector or indicator function
    xcolor,
    mathtools,
    algorithm,
    algpseudocode,
    kbordermatrix,
    authblk,
    enumitem,
    todonotes,
    stmaryrd
}

\usepackage[pdffitwindow=false,
            plainpages=false,
            pdfpagelabels=true,
            pdfpagemode=UseOutlines,
            pdfpagelayout=SinglePage,
            bookmarks=false,
            colorlinks=true,
            hyperfootnotes=false,
            linkcolor=blue,
            urlcolor=blue!30!black,
            citecolor=green!50!black]{hyperref}

\usepackage[bf,normal]{caption}
\graphicspath{{Figures_Revised/}}

\newcommand{\innerprod}[2]{\left\langle #1,\, #2 \right\rangle} % scalar product
\newtheorem{theorem}{Theorem}[section]
\newtheorem{corollary}[theorem]{Corollary}
\newtheorem{lemma}[theorem]{Lemma}
\newtheorem{proposition}[theorem]{Proposition}

\theoremstyle{definition}

\newtheorem{remark}[theorem]{Remark}
\newtheorem{Conjecture}[theorem]{Conjecture}

\makeatletter
\renewcommand*\env@matrix[1][*\c@MaxMatrixCols c]{%
  \hskip -\arraycolsep
  \let\@ifnextchar\new@ifnextchar
  \array{#1}}
\makeatother

\mathtoolsset{centercolon} % for := and =:

\newlength{\continueindent}
\usepackage{etoolbox}
\makeatletter
\newcommand*{\ALG@customparshape}{\parshape 2 \leftmargin \linewidth \dimexpr\ALG@tlm+\continueindent\relax \dimexpr\linewidth+\leftmargin-\ALG@tlm-\continueindent\relax}
\apptocmd{\ALG@beginblock}{\ALG@customparshape}{}{\errmessage{failed to patch}}
\makeatother
\title{Spectral Properties of Effective Dynamics from Conditional Expectations}
\author[1,4]{Feliks N\"uske}
\author[2]{P\'{e}ter Koltai}
\author[3,4]{Lorenzo Boninsegna}
\author[4, 5]{Cecilia Clementi}

\affil[1]{Institute of Mathematics, Universit\"at Paderborn, Germany}
\affil[2]{Department of Mathematics and Computer Science, Freie Universit\"at Berlin, Germany}
\affil[3]{University of California Los Angeles, Institute for Quantitative and Computational Biosciences and Department of Microbiology, Immunology and Molecular Genetics}
\affil[4]{Center for Theoretical Biological Physics and Department of Chemistry, Rice University, USA}
\affil[5]{Department of Physics, Freie Universit\"{a}t Berlin}

\begin{document}
\maketitle
\begin{abstract}
The reduction of high-dimensional systems to effective models on a smaller set of variables is an essential task in many areas of science. For stochastic dynamics governed by diffusion processes, a general procedure to find effective equations is the conditioning approach. In this paper, we are interested in the spectrum of the generator of the resulting effective dynamics, and how it compares to the spectrum of the full generator. We prove a new relative error bound in terms of the eigenfunction approximation error for reversible systems. We also present numerical examples indicating that if Kramers--Moyal (KM) type approximations are used to compute the spectrum of the reduced generator, it seems largely insensitive to the time window used for the KM estimators. We analyze the implications of these observations for systems driven by underdamped Langevin dynamics, and show how meaningful effective dynamics can be defined in this setting.
\end{abstract}

\section{Introduction}
The description of high-dimensional dynamical systems by a reduced set of variables, usually referred to as \textit{coarse graining} or \textit{model reduction}, is of tremendous importance across many different fields of research. Examples range from finance to climate modeling to molecular biology. From the huge body of literature on the subject, we mention in particular the Mori--Zwanzig
formalism~\cite{mori1965transport,Zwa73,ChHaKu00,ChHaKu02,hijon2009mori} as
well as the framework of averaging and homogenization for systems with explicit
multiscale structure~\cite{Pavliotis:2008aa,PAVLIOTIS2007}. Within the field of molecular physics, references \cite{ClementiCOSB,Noe:2017aa,Noid2013,Prinz2011c,Rohrdanz:2011aa,Saunders2013} present important contributions to this line of research. Here, we focus on model reduction for stochastic differential equations (SDEs), and follow another standard approach, which is based on conditioning along level sets of the coarse graining map~\cite{WEINAN2004,Legoll:2010aa,froyland2016trajectory}. For a detailed theoretical analysis of the method in the context of SDEs, please see~\cite{Legoll:2010aa,Zhang:2016aa,Zhang2017,Legoll:2017aa,LELIEVRE2019}.

\noindent For a given coarse grained description of a system, a fundamental question to address is the quality of approximation of the full system by means of the reduced system, as measured by a suitable metric (which usually depends on the problem at hand). In many cases, the approximation of spectral properties of the system's generator is useful in this context. The generator and its associated semigroup, also called Koopman semigroup, are used to describe the time evolution of expectation values of observable functions. For metastable systems, the leading generator eigenpairs provide information on slow modes in the dynamical system. Spectral approximation results for the conditioning approach have been obtained in \cite{Zhang:2016aa,Zhang2017}.

\noindent Another important problem to consider is the analysis and parameter estimation of coarse grained models based on simulation data of the full system. In recent years, a variety of methods has been developed to learn models for the Koopman semigroup off simulation data, see \cite{Schutte:1999aa,Dellnitz:1999aa,Noe2013c,WKR15,MARDT2018,Klus:2018aa,WU020variational} and the references therein. In Ref.~\cite{KNPNCS20}, some of the authors of the present study presented a conceptually simple framework for the data-driven approximation of the Koopman generator, called gEDMD. This framework can also be used to identify and analyze coarse grained models within the context of the conditioning approach. The gEDMD method requires knowledge of the full system parameters. If these parameters are unknown, they need to be replaced by a suitable approximation, such as Kramers--Moyal (KM) formulae. These are based on averages of finite differences at a finite offset (time window). The quality of this approximation as a function of the offset will be addressed in this paper. Even though we focus on KM estimators here, let us mention that a multitude of more advanced methods for parameter estimation of stochastic dynamics are available, please see Ref. \cite{KESSLER2012} for an overview. Spectral methods have been considered in \cite{GOBET2004,Crommelin:2011aa}, while particular attention to the choice of time window has been paid in~\cite{ZHANG2005,PAVLIOTIS2007}. The Kramers--Moyal formulae being among the simplest estimators, we take them as the starting point for our study.

\noindent The third focus of this study is model reduction for systems driven by underdamped Langevin dynamics, which is a widely used model, especially in molecular and biological physics. As the momentum variables of these dynamics often play just an auxiliary role, an interesting question to address is how to define a reduced dynamics that only involves the position state variables. As the conditioning approach does not provide meaningful answers in this case, finding meaningful effective equations remains an open problem in this setting \cite{BITTRACHER2015,Bittracher:2015aa,DUONG2018}.

\noindent In this paper, we report theoretical and numerical results on the issues raised above. The contributions of this study are as follows:

\begin{itemize}
\item Concerning the first problem, we prove a new relative error bound for the approximation of generator eigenvalues by the coarse grained generator, if the dynamics is reversible (Proposition \ref{prop:eigenvalue_error}). This bound shows that a small projection error of the full eigenfunctions with respect to the energy norm is required for a small eigenvalue error. We also derive conditions to ensure that the spectrum of the reduced generator is discrete in the first place (Proposition \ref{prop:discrete_spectrum_lxi}).
\item Concerning the second issue, we present numerical examples indicating that if KM estimators are used within the gEDMD algorithm for reversible systems, on a good set of reaction coordinates, then the resulting eigenvalue estimates seem to be fairly insensitive to the offset used for the KM estimators (Sections \ref{ssec:lemon_slice} and \ref{ssec:toy_molecule}, Conjecture \ref{conj:km_large_s}).
\item Thirdly, we suggest that if the observations of the second part can be confirmed theoretically, it is possible to use KM estimators at \emph{large} offsets to define meaningful effective equations for underdamped dynamics (Corollary \ref{cor:langevinApprox}). The reason is that the statistics of the underdamped process approach those of an overdamped process after a suitable re-scaling of time. We provide successful illustrations of this idea using a toy example and molecular dynamics simulation data of the alanine dipeptide (Sections \ref{ssec:Langevin_toy_model} and \ref{ssec:alanine_dipeptide}).
\end{itemize}

\noindent The rest of this paper is organized as follows: in Section \ref{sec:stochastic_dynamics}, we recap what is needed of the theory of stochastic differential equations, their generators, the conditioning approach, and the data-driven approximation of Koopman generators. In Section \ref{sec:approx_spectrum}, we present and illustrate our spectral approximation result. The technical details of the proofs are deferred to Section~\ref{sec:proofs}. In Section~\ref{sec:km_results}, we present numerical results on the spectrum of gEDMD models based on KM estimators for reversible systems, and conjecture that the observed behaviour can be expected in general. We analyze the implications of this hypothesis for systems driven by underdamped Langevin dynamics in Section~\ref{sec:underdamped}, and provide additional numerical results for this setting. Conclusions and outlook follow in Section~\ref{sec:conclusion}.
 
%%%%%%%%%%%%%%%%%%%%%%%%%%%%%%%%%%%%%%%%%%
\section{Concepts}
\label{sec:stochastic_dynamics}

\subsection{SDEs and Generators}
\label{sec:setting}
In this paper, we consider a reversible Markov process $X_{t}$ attaining values in a domain $\Omega \subset \mathbb{R}^{d}$. The process is governed by the stochastic differential equation
\begin{eqnarray}
dX_{t} &= b(X_{t})dt+ \sigma(X_{t})dB_{t}.\label{eq:ito_sde}
\end{eqnarray}
Here, $B_{t}$ denotes $d$-dimensional Brownian motion, the function $b:\mathbb{R}^{d}\mapsto\mathbb{R}^{d}$ is called the drift, and $\sigma:\mathbb{R}^{d}\rightarrow\mathbb{R}^{d\times d}$ is called the diffusion. We use the notation $a(x) = \sigma(x)\sigma(x)^{T}$ for the covariance matrix of the diffusion. A standard example for dynamics of type Eq.~(\ref{eq:ito_sde}) are the \emph{overdamped Langevin dynamics}
\begin{eqnarray}
dX_{t} = -\frac{1}{\gamma}\nabla V(X_{t})dt+\sqrt{2\beta^{-1}\gamma^{-1}}dB_{t},\label{eq:overdamped_Langevin}
\end{eqnarray}
where $V:\, \Omega \rightarrow \mathbb{R}$ is a scalar function called the \emph{potential}, while $\beta,\,\gamma$ are constants corresponding to the inverse temperature and the friction in physics applications. 

\noindent We assume that $X_{t}$ is ergodic with respect to a unique invariant measure $\mu$ with Boltzmann density $\rho \propto \exp(-F(x))$, where $F$ is called a generalized potential. Also, we assume the diffusion to satisfy a so-called \emph{uniform ellipticity} condition
\begin{equation}
\label{eq:uniform_ellipticity}
0 < \eta_1 \|v\|^2 \leq v^T a(x) v \leq \eta_2 \|v\|^2,
\end{equation}
for constants $\eta_1, \eta_2 > 0$. The invariant measure $\mu$ gives rise to the Hilbert space $L_{\mu}^{2}$ of all square integrable functions with respect to that measure. We can think of these functions as physical observables. The inner product on $L^2_\mu$ is given by
\begin{equation}
\label{eq:inner_product_mu}
\langle \psi, \tilde{\psi} \rangle_{\mu} = \int_{\Omega}\psi(x) \tilde{\psi}(x)\,\mathrm{d}\mu(x) = \int_\Omega \psi(x) \tilde{\psi}(x) \rho(x) \,\mathrm{d}x.
\end{equation}
For a fixed time window $t \geq 0$, and an observable function $\psi \in L^2_\mu$, the \emph{Koopman operator} $\mathcal{K}_{t}$ describes the evolution of the expectation value of $\psi$ by means of the dynamics \eqref{eq:ito_sde}:
\begin{equation*}
\mathcal{K}_{t} \psi(x) = \mathbb{E}^{x}\left[\psi(X_{t})\right],
\end{equation*}
where $\mathbb{E}^x[\cdot]$ denotes expectation given that the dynamics starts deterministically at $x$. The \emph{infinitesimal generator} $\mathcal{L}$ of the Markov process $X_t$ is then defined as a formal time-derivative of this expectation value:
\begin{equation}
\label{eq:definition_generator}
\mathcal{L}\psi(x) = \frac{d}{dt}\mathbb{E}^x\left[ \psi(X_t) \right]\vert_{t=0}.
\end{equation}
It follows that, by describing the system in terms of the expectations $\mathcal{K}_t \psi$, the non-linear dynamical system \eqref{eq:ito_sde} turns into a linear, but infinite-dimensional system with differential equation
\begin{equation*}
\frac{d}{dt} \mathcal{K}_t \psi = \mathcal{L}\mathcal{K}_t \psi.
\end{equation*}
For this reason, the Koopman operators $\mathcal{K}_t$ and their generator $\mathcal{L}$ have been studied extensively in past decades. A technical subtlety arising from this infinite-dimensional description is that the time-derivative \eqref{eq:definition_generator} is not well-defined for \emph{all} $\psi \in L^2_\mu$. For smooth and compactly supported functions however, stochastic calculus shows that $\mathcal{L}$ is well-defined, and acts as a second order differential operator:
\begin{align}
\mathcal{L}\psi(x) &= \sum_{i=1}^d b_i(x) \frac{\partial \psi(x)}{\partial x_i} + \frac{1}{2}\sum_{i,j=1}^d a_{ij}(x) \frac{\partial^2 \psi(x)}{\partial x_i \partial x_j} \\
&= b(x)\cdot\nabla \psi(x)+\frac{1}{2} a(x):\nabla^{2} \psi(x).\label{eq:generator_differential_op}
\end{align}
Here, $\nabla^{2}\psi$ is the Hessian matrix of the function $\psi$, and the colon $:$ denotes the Frobenius inner product between matrices, i.e., $A:B = \sum_{i,j} A_{ij}B_{ij}$. For the same class of functions, the generator is symmetric with respect to the inner product \eqref{eq:inner_product_mu}, and satisfies the important equality
\begin{eqnarray}
\langle\mathcal{L}\psi, \tilde{\psi} \rangle_{\mu} & = & -\frac{1}{2}\int_{\Omega}a(x)\nabla \psi(x)\cdot\nabla \tilde{\psi}(x)\,\mathrm{d}\mu(x),\label{eq:scalar_prod_generator}
\end{eqnarray}
which requires only first order derivatives. The negative of the right-hand side of Eq. \eqref{eq:scalar_prod_generator} is called the \emph{quadratic form}
\begin{equation}
\label{eq:def_quadratic_form}
\mathcal{Q}(\psi, \tilde{\psi}) := \frac{1}{2}\int_{\Omega} a(x)\nabla \psi(x)\cdot\nabla \tilde{\psi}(x)\,\mathrm{d}\mu(x).
\end{equation}
Tools from functional analysis \cite{BAKRY2013,Davies:1982aa} can be used to define \eqref{eq:def_quadratic_form} on a larger set of functions $\mathbb{V}_\mu$, usually called the \emph{form domain}. Because of \eqref{eq:uniform_ellipticity}, $Q$ defines an inner product on the form domain, and $\mathbb{V}_\mu$ in fact turns into a Hilbert space with \emph{energy norm}
\begin{equation*}
\|\cdot \|_\mathcal{Q} = \mathcal{Q}(\cdot, \cdot)^{1/2}.
\end{equation*}
The reason to introduce all these concepts is that the energy norm will serve as error measure for the main results of this study.

\subsection{Spectral Decomposition}
We are particularly interested in eigenvalues and eigenfunctions of the (negative) generator $-\mathcal{L}$. Because of \eqref{eq:scalar_prod_generator}, the spectrum of $-\mathcal{L}$ must be part of the non-negative real axis. We will further assume that there is a complete set of eigenfunctions (i.e., they form a basis of $L_{\mu}^{2}$) corresponding to discrete eigenvalues. In other words, there are functions $\psi_0, \psi_{1},\psi_{2},\ldots$ and non-negative numbers $\kappa_{0} < \kappa_1 < \kappa_2 < \ldots$ such that
\begin{eqnarray}
-\mathcal{L}\psi_{i} = \kappa_{i}\psi_{i}.\label{eq:eigenfunction_generator}
\end{eqnarray}
Conditions for the existence of a completely discrete spectrum are discussed in section \ref{sec:proofs}. The assumption that all eigenvalues $\kappa_i$ are distinct is for simplicity only, especially with regards to the analysis in section~\ref{sec:proofs}. It follows again from \eqref{eq:scalar_prod_generator} that $\kappa_0 = 0$, and $\psi_0 \equiv 1$ is the constant function.

\noindent The physical significance of these eigenpairs is that by \cite[Ch 2, Thm 2.4.]{Pazy:1983aa}, the eigenfunctions $\psi_{i}$ are also eigenfunctions of the Koopman operators $\mathcal{K}_{\tau}$ for all $t \geq 0$, corresponding to eigenvalues
\begin{equation}
\label{eq:eigenvalue_koopman_op}
\lambda_{i}(t) = e^{-\kappa_{i}t}.
\end{equation}
Due to the exponential decay of all $\lambda_{i}(t)$, it is common to refer to the $\kappa_{i}$ as rates, and to their reciprocals as \emph{implied timescales}
\begin{eqnarray}
t_{i} = \frac{1}{\kappa_{i}}.\label{eq:def_its}
\end{eqnarray}
In many applications, including molecular dynamics, we expect to find a number $K$ of dominant rates $0 < \kappa_{1}< \ldots <\kappa_{K}\ll\kappa_{K+1}$ separated from all others. These dominant spectral components are of particular interest as they are related to metastability, that is, the existence of long-lived macrostates such that transitions between those states are rare events \cite{Davies1982a,Dellnitz:1999aa,Deuflhard2000f}.

\subsection{Dimensionality Reduction}
\label{sec:projected_dynamics}
The main topic of this study is the effect of dimensionality reduction on the generator eigenvalues $\kappa_i$ introduced above. Following the notation of Refs. \cite{Legoll:2010aa,Zhang:2016aa,Zhang2017}, we consider a smooth coarse graining function $\xi$, which maps the state space $\Omega \subset \mathbb{R}^{d}$ onto a lower-dimensional space $\hat{\Omega}\subset \mathbb{R}^{m}$, where $m\leq d$. For a position $z \in \hat{\Omega}$ in reduced space, we denote the marginal probability distribution of the invariant measure $\mu$ by $\nu$, and assume it possesses a corresponding density function $\vartheta(z)$. The \emph{conditional expectation operator}
\begin{equation*}
\mathcal{P}\psi(z) = \mathbb{E}^{\mu}\left[ \psi(x) \mid \xi(x) = z \right]
\end{equation*}
computes the stationary average of a function $\psi$ defined on $\Omega$, conditional to $\xi$ attaining a fixed value $z \in \hat{\Omega}$. Using the Dirac $\delta$-function, we can informally write the above expression as
\begin{equation*}
\mathcal{P}\psi(z) = \frac{1}{\vartheta(z)} \int_\Omega \psi(x) \delta(\xi(x) - z) \, \mathrm{d}\mu(x).
\end{equation*}

\noindent Consider the space $L^2_\nu$ of physical observables on reduced space $\hat{\Omega}$:
\begin{equation*}
L^2_\nu = \{\varphi : \hat{\Omega} \rightarrow \mathbb{R}, \, \int_{\hat{\Omega}} \varphi^2(z) \,\mathrm{d}\nu(z) < \infty  \}.
\end{equation*}
By the concatenation $\varphi \circ \xi$, every such function can be viewed as a function on full state space $\Omega$. In fact, $L^2_\nu$ can be exactly identified as the subspace of functions in $L^2_\mu$ which depend only on the value of $z = \xi(x)$, with the conditional expectation operator acting as orthogonal projection onto this subspace~\cite{Zhang:2016aa}. Note that unless $\xi$ is constant, $L_{\nu}^{2}$ is an infinite-dimensional subspace.

\noindent This motivates consideration of the projected generator
\begin{eqnarray}
\mathcal{L}^{\xi} &= \mathcal{PLP}.\label{eq:def_projected_generator}
\end{eqnarray}
As discussed comprehensively in \cite{Zhang:2016aa}, this operator retains the shape of the generator of a reversible Markov process $Z_t$ on $\hat{\Omega}$, as for smooth compactly supported functions $\varphi \in L^2_\nu$, we have that
\begin{equation}
\label{eq:proj_generator_diff_op}
\mathcal{L}^{\xi}\varphi = \mathcal{P}\left[\mathcal{L}\xi\right]\cdot\nabla_{z}\varphi +\frac{1}{2}\mathcal{P}\left[\nabla\xi^{T}a\nabla\xi\right]:\nabla_{z}^{2}\varphi.
\end{equation}
In the above equation, $\mathcal{L}\xi$ is an $m$-dimensional vector, each entry containing the application of $\mathcal{L}$ to each component of $\xi$, and $\nabla\xi$ is the $d \times m$ Jacobian matrix of $\xi$. The coefficients
\begin{align}
b^{\xi}(z) &= \mathcal{P}(\mathcal{L}\xi)(z),&
a^{\xi}(z) &= \mathcal{P}\left(\nabla\xi^{T} a \nabla\xi \right)(z)\label{eq:projected_drift_diff}
\end{align}
serve as \emph{effective drift} and \emph{effective diffusion} for the process $Z_t$, respectively. It can also be shown that $Z_t$ is ergodic with respect to $\nu$ \cite{Zhang:2016aa}, and we can associate to $\mathcal{L}^\xi$ a form domain $\mathbb{V}_\nu$ with an effective quadratic form
\begin{equation}
\label{eq:effective_quad_form}
\mathcal{Q}^\xi(\varphi, \tilde{\varphi}) = \frac{1}{2}\int_{\hat{\Omega}} a^\xi(z) \nabla_z \varphi(z) \cdot \nabla_z \tilde{\varphi}(z) \,\mathrm{d}\nu(z).
\end{equation}

\noindent If the projected generator also possesses a discrete spectrum with eigenvalues $\omega_i, \, i = 0, 1, \ldots$, comparison of those eigenvalues with the original ones $\kappa_i$ provides information about how well the effective dynamics $Z_t$ retain the relaxation processes of the original process. Our results on this topic are presented in sections \ref{sec:approx_spectrum} and~\ref{sec:proofs}.

\subsection{Galerkin Approximation}
\label{sec:MSM_VAC}
The numerical approximation of the eigenfunctions $\psi_i, \, i=1,\ldots, K$ is often achieved by Galerkin projection, i.e., orthogonal linear projection of the generator to a finite-dimensional subspace. After choosing such a space $\mathbb{W} \subset \mathbb{V}_\mu$, with a basis set $\{\phi_i\}_{i=1}^N$, the Galerkin approach consists of finding $\hat{\psi}_i \in \mathbb{W}$ such that
\begin{align}
\label{eq:galerkin_ev_problem}
-\innerprod{\mathcal{L}\hat{\psi}_i}{\phi_j}_\mu = Q(\hat{\psi}_i, \phi_j) &= \hat{\omega}_i \innerprod{\hat{\psi}_i}{\phi_j}_{\mu}\quad \forall 1 \leq j \leq N.
\end{align}
Eq. \eqref{eq:galerkin_ev_problem} is a generalized matrix eigenvalue problem. The approximate eigenvalues $\hat{\omega}_i$ are also called \emph{Ritz values} associated to~$\mathbb{W}$. If $\mathbb{W}$ is chosen as a subspace of $\mathbb{V}_\nu$ (given a set of reduced variables as described in the previous section), it was shown in \cite{Zhang:2016aa} that \eqref{eq:galerkin_ev_problem} serves both as a weak form for $\mathcal{L}$ in $\mathbb{V}_\mu$ and for $\mathcal{L}^\xi$ in $\mathbb{V}_\nu$. Importantly, the min-max-principle implies that for any such subspace
\begin{equation}
\label{eq:min_max_L_Lxi}
\kappa_i \leq \omega_i \leq \hat{\omega}_i.
\end{equation}

\noindent Moreover, if simulation data $\{x_m\}_{m=1}^M$ of the full process \eqref{eq:ito_sde}, approximately sampling the invariant measure $\mu$, is available, the following empirical estimators
\begin{align}
\label{eq:gedmd_estimators}
\mathcal{Q}(\phi_i, \phi_j) &\approx \frac{1}{2M}\sum_{m=1}^M \nabla \phi_i(x_m)^T a(x_m) \nabla \phi_j(x_m), & \innerprod{\phi_i}{\phi_j}_\mu &\approx \frac{1}{M}\sum_{m=1}^M \phi_i(x_m)\phi_j(x_m),
\end{align}
will converge to the terms in Eq. \eqref{eq:galerkin_ev_problem} in the limit of infinite data. This method has been called \emph{generator Extended Dynamic Mode Decomposition} (gEDMD) \cite{KNPNCS20}. Consequently, for a subspace $\mathbb{W} \subset \mathbb{V}_\nu$ comprised of functions on reduced space, gEDMD will simultaneously approximate the eigenvalues of $\mathcal{L}$ and $\mathcal{L}^\xi$ by \eqref{eq:min_max_L_Lxi}. Note that in this setting, simulation data of the full dynamics \eqref{eq:ito_sde} can still be used in Eq. \eqref{eq:gedmd_estimators}, it is not necessary to simulate the effective dynamics defined by $\mathcal{L}^\xi$ first. As the gradients in \eqref{eq:gedmd_estimators} are taken with respect to the full state variables $x$, we have to apply the chain rule when evaluating \eqref{eq:gedmd_estimators} for basis functions $\phi_i \in \mathbb{W} \subset \mathbb{V}_\nu$: $\nabla \phi_i(\xi(x)) = \nabla_z \phi_i(z)^T \nabla \xi^T(x)$.

\noindent Note that the estimator for $\mathcal{Q}$ in \eqref{eq:gedmd_estimators} exploits reversibility of the process. An alternative estimator, which can also be applied to non-reversible processes, is
\begin{align}
\label{eq:gedmd_estimators_non_rev}
-\innerprod{\mathcal{L}\phi_i}{\phi_j}_\mu &\approx -\frac{1}{M}\sum_{m=1}^M \left[b(x_m) \nabla \phi_i(x_m) + \frac{1}{2} a(x_m) : \nabla^2 \phi_i(x_m)\right] \phi_j(x_m).
\end{align}
We will require this last equation in the next section.

\subsection{Kramers--Moyal Estimators}
\label{ssec:km_introduction}
The parameters $b^\xi$ and $a^\xi$ in \eqref{eq:projected_drift_diff} involve integrals over non-linear manifolds in high-dimensional space, and they are rarely used in practice for this reason. For a process $\{X_t\}_{t \geq 0}$ and a positive \emph{offset} $s > 0$, define the projected first order finite difference at time $t$ by 
\begin{equation*}
d^\xi_s(X_t) = \xi(X_{t+s}) - \xi(X_t).
\end{equation*}
With this notation, the basic \emph{Kramers--Moyal (KM) formulae} \cite{Risken1989} for the approximation of $b^\xi$ and $a^\xi$ are given by:
\begin{align}
\label{eq:Kramers_Moyal}
b^{\xi}(z) & = \lim_{s\rightarrow0}\mathbb{E}^\mu \left[\frac{1}{s} d^\xi_s(X_0) \Big\vert \xi(X_{0}) = z\right], &
a^{\xi}(z) &= \lim_{s\rightarrow0}\mathbb{E}^\mu \left[\frac{1}{2s} d^\xi_s(X_0) \otimes d^\xi_s(X_0)  \Big\vert \xi(X_{0}) = z \right].
\end{align}
We note that many more sophisticated approximations can be found in the literature.

\noindent KM estimators can also be used in the context of gEDMD.  If the parameters $b$ and $a$ of the original process \eqref{eq:ito_sde} are unknown, gEDMD can be re-formulated upon replacing $b$ and $a$ by the first order finite differences at all data points in \eqref{eq:gedmd_estimators_non_rev}. If $s > 0$ corresponds to a multiple of the integration time step in a discrete trajectory $\{x_m\}_{m=1}^M$, then \eqref{eq:gedmd_estimators_non_rev} can be converted to:

\begin{equation}
\label{eq:gedmd_km}
-\innerprod{\mathcal{L}\phi_i}{\phi_j}_\mu \approx -\frac{1}{M}\sum_{m=1}^M \left[\frac{1}{s} d^\xi_s(x_m)  \nabla_z \phi_i(x_m) + \frac{1}{2s} (d^\xi_s(x_m) \otimes d^\xi_s(x_m)) : \nabla^2_z \phi_i(x_m)\right] \phi_j(x_m).
\end{equation}

\noindent This estimator is consistent for $s \rightarrow 0$ and $M \rightarrow \infty$. There is a dual meaning to eigenpairs derived from this approximation: on one hand, they serve as an approximation to the eigenpairs of $\mathcal{L}^\xi$, at least for small offsets $s$. On the other hand, these are also approximations to the spectrum of the generator of a non-reversible coarse grained SDE with drift and diffusion given by \eqref{eq:Kramers_Moyal}, see again \cite{KNPNCS20} for a detailed exposition. In Secs. \ref{sec:km_results} and \ref{sec:underdamped}, we present a numerical study on the effect of the offset $s$ on the spectrum of this generator. We note that similar approximations can certainly be built on more advanced estimators than \eqref{eq:Kramers_Moyal}, but the KM formulae will suffice for this study.

\section{Spectral Properties of the Projected Generator}
\label{sec:approx_spectrum}

\subsection{Summary of spectral properties}
\label{ssec:summary}
The first major result of this study concerns the approximation error for the dominant eigenvalues $\kappa_{1}, \kappa_2, \ldots$ by corresponding eigenvalues of $\mathcal{L}^\xi$. We only provide a high-level summary of these results here, while the technically more involved statements and their proofs can be found in section~\ref{sec:proofs}.

\noindent First, we introduce conditions to ensure the spectrum of the effective generator $\mathcal{L}^\xi$ is also discrete, see Proposition~\ref{prop:discrete_spectrum_lxi}. We then show in Proposition~\ref{prop:eigenvalue_error} and Corollary~\ref{cor:simpler_bound} that the \emph{relative} eigenvalue error 
\begin{equation*}
E_i = \frac{\omega_{i}-\kappa_{i}}{\omega_{i}}
\end{equation*}
can be bounded in terms of the energy norm of the projection residual $\mathcal{P}_Q^{\perp}\psi_i = (\mathcal{I} - \mathcal{P}_Q)\psi_i$ of the corresponding eigenfunctions, with $\mathcal{P}_Q$ denoting the $\mathcal{Q}$-orthogonal projection onto $\mathbb{V}_\nu$:
\begin{equation}
\label{eq:error_bound_summary}
E_i = \frac{\omega_{i}-\kappa_{i}}{\omega_{i}} \leq C \|\mathcal{P}_\mathcal{Q}^\perp \psi_i \|^2_\mathcal{Q}.
\end{equation}
In other words, if the eigenfunction $\psi_i$ can be written as a function of the reduced variables $z$, up to a small error, then we can expect the eigenvalue $\kappa_i$ to be reproduced well by the effective dynamics on $z$. However, as the error is measured by the energy norm, Proposition \ref{prop:eigenvalue_error} shows that not only the eigenfunction $\psi_i$, but also its first order derivatives must be approximated well by functions of $z$ alone. In the next section \ref{ssec:illustrate_bound}, we show that this is not merely an academic condition, but indeed necessary.

\noindent Our result improves on existing ones in two ways. First, in Ref.~\cite[Theorem~2]{Zhang2017}, it was shown that the \emph{absolute} eigenvalue
error of the projected generator is small if~$\smash{ \|\mathcal{LP}^{\perp}\psi_{i}\|_{L_{\mu}^{2}} }$ and~$\smash{ \|\mathcal{P}^{\perp}\psi_{i}\|_{L_{\mu}^{2}} }$ are small. Proposition~\ref{prop:eigenvalue_error} (Corollary~\ref{cor:simpler_bound}) complements these results in the sense that it bounds the \emph{relative} error of eigenvalues (timescales), which is a more practical error measure for eigenvalues close to zero, i.e., large timescales. For a more detailed elaboration on the relationship of the projection error and the relative error of timescales, please refer to the text after Corollary~\ref{cor:simpler_bound}.
Second, our bound~(\ref{eq:error_bound_rates}) is less restrictive than the conditions in~\cite{Zhang2017}, as it uses the energy norm involving only first order derivatives, while the term $\smash{ \|\mathcal{LP}^{\perp}\psi_{i}\|_{L_{\mu}^{2}} }$ necessarily requires second derivatives. In fact, the bound assumed in~\cite[Theorem~2]{Zhang2017} implies our bound up to another multiplicative constant (thus, it is more restrictive), see Lemma~\ref{lem:statement_adn_lemma}.

\subsection{Illustration of the Error Bound}
\label{ssec:illustrate_bound}
\begin{figure}[htb]
\begin{centering}
\includegraphics[width=.9\textwidth]{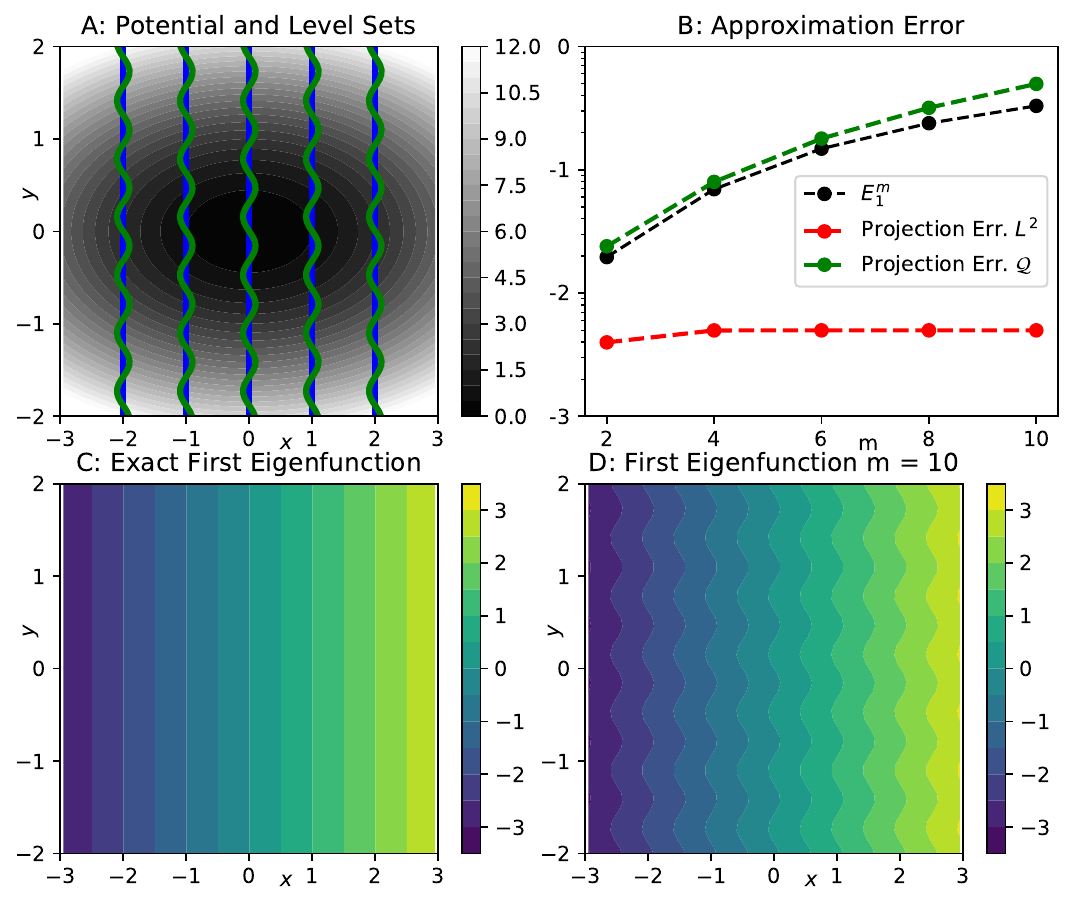}
\par\end{centering}
\centering{}\caption{Illustration of Proposition \ref{prop:eigenvalue_error} by means of a two-dimensional Ornstein-Uhlenbeck process, and one-dimensional reaction coordinates $\xi^m(x, y) = x + 0.1\sin(my)$. A: Selected level sets $\xi^{-1}(z)$ for $m = 0$  (blue) and $m = 10$ (green), with a contour of the potential in the background. B: Red: $L^2_\mu$-error between exact first eigenfunction $\psi_{1,0} = x$, and the approximate slowest eigenfunction $\hat{\psi}^m_1$, computed by Galerkin projection onto the space of the first ten Hermite polynomials $\psi_i(z)$, where $z$ is the reaction coordinate $\xi^m$. Green: the same error, but measured using the energy norm. Black: Relative eigenvalue error $E^m_1$ corresponding to the same approximation. The vertical axis is labeled by the decadic logarithm. C: Contour of the exact slowest eigenfunction $\psi_{1, 0} = x$. D: Contour of the approximate slowest eigenfunction $\hat{\psi}^m_1$ for $m = 10$.\label{fig:illustration_error_bound}}
\end{figure}

\noindent In order to provide an illustration of the error bound, we consider a two-dimensional Ornstein--Uhlenbeck process, that is, overdamped Langevin dynamics \eqref{eq:overdamped_Langevin} with potential (see Fig.~\ref{fig:illustration_error_bound} A):
\begin{eqnarray}
V(x,y) & = \frac{1}{2}(\alpha_x x^2 + \alpha_y y^2).\label{eq:2d_OU_potential}
\end{eqnarray}

\noindent We set $\alpha_x = 1, \alpha_y = 5, \beta = \gamma = 1.0$. The eigenvalues and eigenfunctions of this system are known analytically. For all integers $r, s \geq 0$, we have eigenvalues $\kappa_{r,s} = r\alpha_x + s\alpha_y$, with eigenfunctions $\psi_{r,s} = \frac{1}{\sqrt{r!s!}} H_r(x) H_s(y)$, using one-dimensional Hermite polynomials $H_i$. The first four non-zero eigenvalues correspond to eigenfunctions which are constant in $y$-direction. The first of these eigenfunctions, $\psi_{1,0} = x$, is shown in Fig. \ref{fig:illustration_error_bound} C. We now consider a family of reaction coordinates
\begin{equation*}
\xi^m(x, y) = x + 0.1\sin(m y).
\end{equation*}
For $m = 0$, the reaction coordinate $\xi^0(x, y) = x$ perfectly captures the first eigenfunction $\psi_{1,0}$. For positive $m$, however, the level sets  $\xi^{-1}(z)$ of the coarse graining map oscillate within a vertical strip of width 0.1 around $x = z$, see Fig. \ref{fig:illustration_error_bound} A for a comparison of the level sets at $m = 0$ and $m = 10$. Due to these small scale oscillations, we still expect to find a relatively small projection error of the eigenfunction $\psi_{1,0}$, if measured by the $L^2_\mu$-norm, but an increasingly larger error if the energy norm (involving derivatives) is employed.

\noindent In order to estimate these projection errors, we use finite-dimensional subspaces $\mathbb{W}^m \subset \mathbb{V}^m_{\nu}$, where $\mathbb{V}^m_{\nu}$ is the form domain corresponding to reaction coordinate $\xi^m$. The subspaces are spanned by the first ten Hermite polynomials $H_i(z)$, which exactly capture the slowest eigenfunction for $m = 0$. For each subspace, we calculate the Galerkin matrices in Eq. \eqref{eq:galerkin_ev_problem} by numerical integration, and then extract the first non-trivial eigenvalues $\hat{\omega}^m_{1}$ and eigenfunctions $\hat{\psi}^m_{1}$. We calculate the relative errors $E^m_{1} = \frac{\hat{\omega}^m_{1} - \kappa_{1,0}}{\hat{\omega}^m_{1}}$ and the eigenfunction approximation errors $\|\psi_{1, 0} - \hat{\psi}^m_{1}\|^2$, measured using both the norms on $L^2_\mu$ and the energy norm. The results are shown in Fig. \ref{fig:illustration_error_bound} B. We observe that the $L^2_\mu$-error remains almost constant as $m$ increases, while the energy norm error increases steadily, reflecting the increasingly oscillating shape of the approximate eigenfunctions $\hat{\psi}^m_1$ (see Fig. \ref{fig:illustration_error_bound} D). In agreement with our error estimate, the energy norm approximation error provides a fairly tight bound for the relative eigenvalue error $E^m_1$. It should be noted that the quantities shown here only provide upper bounds for the relative error $E_1$ and the projection error $\delta_1^2$, but they suffice for the purpose of illustration.

\noindent In summary, this example confirms the eigenvalue error bound provided in Proposition \ref{prop:eigenvalue_error}, and it also highlights the importance of capturing the derivatives of generator eigenfunctions of interest when selecting a reaction coordinate for coarse graining.

\section{Spectral Properties and Kramers--Moyal Estimators}
\label{sec:km_results}
The second part of this work is a numerical study on the effect of employing Kramers--Moyal type approximations when estimating spectral properties of projected generators. As explained in section \ref{ssec:km_introduction}, KM estimators can be incorporated into the gEDMD method by means of Eq. \eqref{eq:gedmd_km} if the full system parameters $b$ and $a$ are not known. Again, we stress that using gEDMD to calculate eigenvalues in this way possesses a dual meaning: it is an approximation to the eigenvalues of $\mathcal{L}^\xi$, but also an approximation to the eigenvalues of a non-reversible dynamics with parameters given by \eqref{eq:Kramers_Moyal}. We find that for projections onto a good set of reaction coordinates $\xi$, i.e., coordinates which are known to capture the slow dynamics of a system well, it is possible to use a surprisingly large offset parameter $s$ in Eq. \eqref{eq:gedmd_km}.

\subsection{Methods}
\label{ssec:methods}
In all of the following examples, we generate long realizations of the system under investigation, by employing the Euler-Maruyama method with discrete integration time step $\Delta_t$, for $M$ steps, such that the total simulation time equals $M\Delta_t$. The only exception is the molecular example in section \ref{ssec:alanine_dipeptide}, where a molecular dynamics code was used to generate the data, see the references given there. For the application of gEDMD on one-dimensional domains, we either use Gaussian basis functions $\phi_i$, or periodic Gaussians $\phi^p_i$ if the domain is periodic:
\begin{align}
\label{eq:gaussian_fct}
\phi_i(z) &= \exp\left[-\frac{1}{2 \rho}(z - z_i)^2\right], & \phi^p_i(z) &= \exp\left[-\frac{1}{2 \rho}\sin^2(\frac{1}{2}(z - z_i))\right],
\end{align}
with bandwidth $\rho$ and centers $z_i \in \hat{\Omega}$. On two-dimensional domains, we use products of the univariate functions defined above, centered on a regular grid.

\noindent With these basis functions, the Galerkin matrices for gEDMD are calculated in different ways. As a reference, we use estimators \eqref{eq:gedmd_estimators}, which require knowledge of the full system parameters. In addition, we use the estimator \eqref{eq:gedmd_km} with a series of offsets $s > 0 $. We then solve the generalized eigenvalue problem \eqref{eq:galerkin_ev_problem} for each of these cases. Eigenvalue estimates thus obtained are denoted by $\hat{\omega}_i^0$ and $\hat{\omega}_i^s$, respectively. We keep track of the relative error
\begin{equation}
\label{eq:relative_err_ev}
E_i^s = \frac{|\hat{\omega}^s_i - \hat{\omega}^0_i|}{\hat{\omega}^s_i},
\end{equation}
and also monitor the reciprocals $\hat{t}_i^0 = \left(\hat{\omega}_i^0\right)^{-1}, \,\hat{t}_i^s = \left(\hat{\omega}_i^s\right)^{-1}$, which serve as estimates of the implied timescales \eqref{eq:def_its}. Additionally, we also extract estimates of the $K$ dominant eigenfunctions from each of these eigenvalue problems, and apply the PCCA method \cite{Deuflhard:2005aa} to determine metastable decompositions of the domain based on these eigenfunctions. PCCA returns $K$ membership functions $\chi_j(z), j = 1,\ldots, K$, such that $\sum_{j = 1}^K \chi_j(z) = 1$ at all points $z$, and $\chi_j(z)$ indicates the degree of membership of each point $z$ to metastable state $j$.

\subsection{Lemon Slice Potential}
\label{ssec:lemon_slice}
\begin{figure}
\begin{centering}
\includegraphics{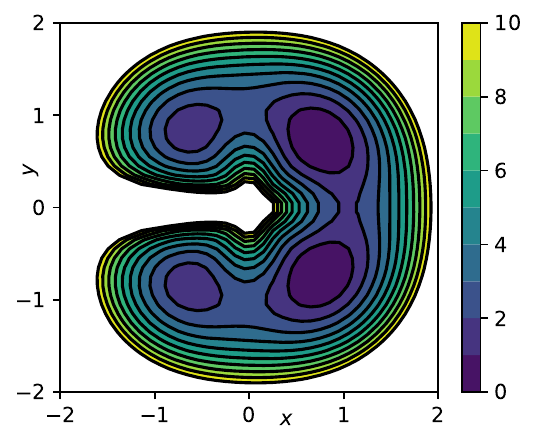}
\par\end{centering}
\caption{Contour plot of the lemon slice potential Eq. (\ref{eq:lemon_slice}).\label{fig:lemon_slice_potential}}
\end{figure}

\begin{figure}
\begin{centering}
\includegraphics[width=.95\textwidth]{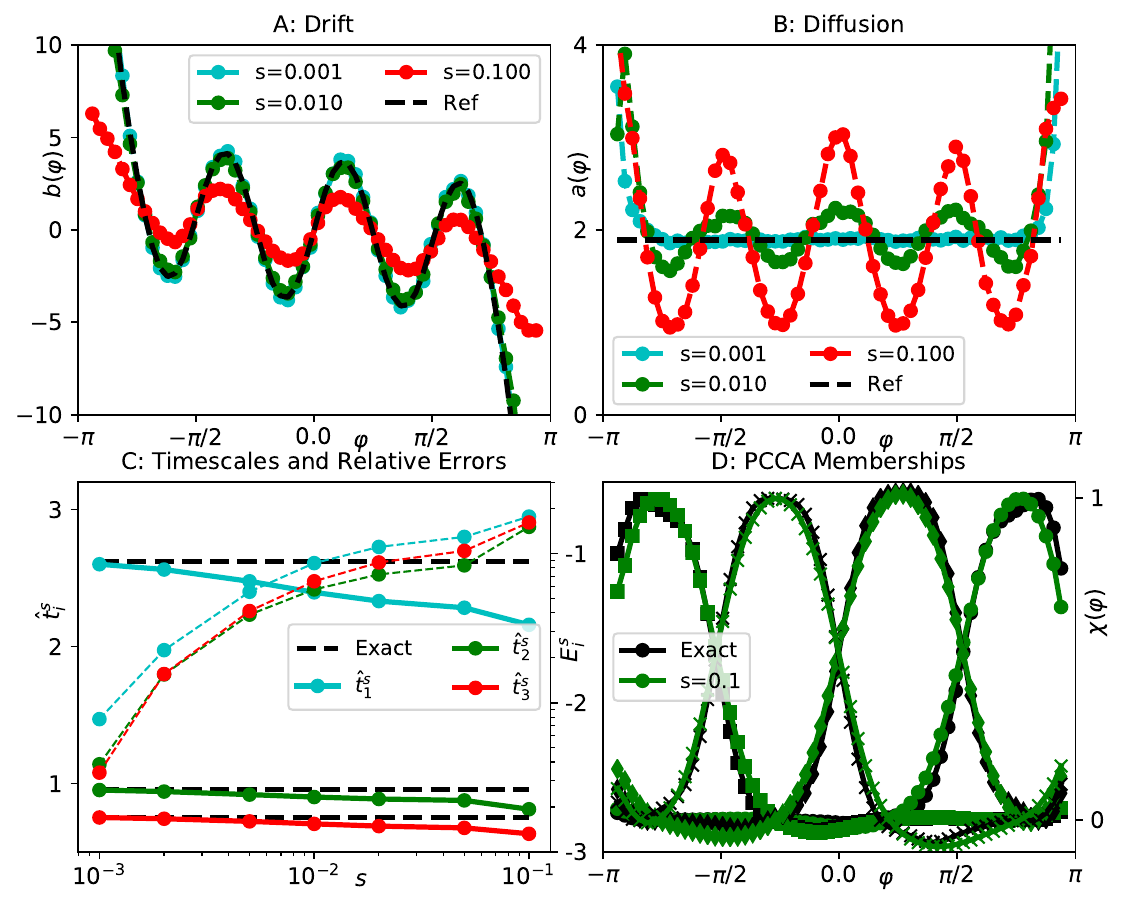}
\par\end{centering}
\centering{}\caption{Analysis of effective dynamics on the polar angle for the lemon slice potential. A: Numerical estimates of effective drift for different values of the offset $s$, compared to the reference in black. B: The same for the effective diffusion. C: Implied timescales $\hat{t}_{1}^s,\, \hat{t}_{2}^s,\, \hat{t}_{3}^s$ extracted from gEDMD models using KM formulae at various offsets $s$ (solid lines), compared to the results of applying gEDMD with full system parameters (dashed black lines). We also show the relative errors $E_i^s$ \eqref{eq:relative_err_ev} for all three timescales (thin dashed lines, scale on the right, labeled by decadic logarithm). D) Four metastable membership functions generated by the PCCA method, extracted from a gEDMD model at offset $s = 0.1$ (green) and using exact system parameters (black).\label{fig:lemon_slice}}
\end{figure}

We consider overdamped Langevin dynamics Eq.~(\ref{eq:overdamped_Langevin}) in the ``lemon slice'' potential
\begin{equation}
V(r,\varphi) = \cos(4 \varphi) + \frac{1}{\cos(0.5\varphi)} + 10(r-1)^{2}+\frac{1}{r},\label{eq:lemon_slice}
\end{equation}
where $r,\varphi$ are two-dimensional polar coordinates, at inverse temperature $\beta=1.0$ and friction $\gamma = 1.0$. A contour of the potential is shown in Fig.~\ref{fig:lemon_slice_potential}. Note that the second and the last term in \eqref{eq:lemon_slice} impose an infinite barrier along the negative $x$-axis and at the origin. The slow dynamics of this system correspond to transitions between the four main minima of the potential $V$, we therefore find three dominant timescales $t_{1} \approx 2.6,\, t_2 \approx 0.95,\, t_{3} \approx 0.75$, see \cite{KNPNCS20} for a previous analysis of the same example. Hence, we can select the polar angle $\varphi$ as a suitable reaction coordinate $\xi(x,y) =  \varphi(x,y)$. In this case, the effective drift $b^{\xi}$ and diffusion $a^{\xi}$ can even be calculated analytically, see appendix~\ref{sec:appendix_lemon_slice}, they are indicated by the black lines in Figs. \ref{fig:lemon_slice} A and B. The data set we use comprises $M = 5\cdot 10^6$ data points at integration time step $\Delta_t = 10^{-3}$.

\noindent Since applying gEDMD with a positive offset in \eqref{eq:gedmd_km} corresponds to approximating the generator of an SDE with coefficients \eqref{eq:Kramers_Moyal}, and we also have analytical expressions for the exact effective parameters $b^\xi$ and $a^\xi$, we first provide a comparison between these analytical parameters and histogrammed estimates of Eq. \eqref{eq:Kramers_Moyal}. We find in Figs. \ref{fig:lemon_slice} A and B, that the exact drift and diffusion are recovered well for small $s$, as expected, while very different results are obtained as $s$ increases. In particular, the effective diffusion is no longer constant as a function of $\varphi$ for large $s$. 

\noindent Next, we compare the dominant spectra of the generators corresponding to these different dynamics. To this end, we employ gEDMD with fifteen Gaussian basis functions \eqref{eq:gaussian_fct}, centered at equal distance between $z_i = -2.8$ and $z_i = 2.8$, each of width $\rho = 0.1$, as described in section \ref{ssec:methods}. In Fig.~\ref{fig:lemon_slice} C, we show the first three timescales $\hat{t}_i^s$ compared to their reference values, and also the relative errors $E_i^s$, see \eqref{eq:relative_err_ev}. As $s$ increases by about two orders of magnitude, the relative errors remain within a ten to twenty percent margin around the references, which is generally acceptable from a practical point of view. As shown in Fig.~\ref{fig:lemon_slice} D, the membership functions $\chi_j, j=1, \ldots, 4$ generated by the PCCA method for $s = 0.1$, barely differ from the ones computed using full system parameters. We therefore conclude that the dominant spectrum is approximately retained by all models up to $s = 0.1$. In other words, the structure of the spectrum is unchanged as long as the offset is clearly smaller than the slow timescales, even though the KM estimates for drift and diffusion are very different from the analytical coefficients $b^\xi$ and $a^\xi$.

\subsection{Prototypical Molecular Potential}
\label{ssec:toy_molecule}
In order to confirm the observations made in the previous section, we study a more complex example. The system is designed to mimic a small molecule consisting of five atoms. The three-dimensional Euclidean coordinates of all five atoms thus define the system's fifteen-dimensional state space. We imagine these five atoms to be linked by bonds like a chain. The dynamical model is again the overdamped Langevin dynamics \eqref{eq:overdamped_Langevin}, with a potential energy comprised of typical molecular interactions, namely harmonic bond, bond angle, and dihedral potentials. The precise parameters can be found in appendix \ref{sec:appendix_molecule}. Note, however, that the system is only qualitatively similar to a small molecule, since parameter values and the units of time and energy do not correspond to physical values. We also do not include any solvent molecules or velocities. The data set comprises $M = 2\cdot 10^6$ data points at integration time step $\Delta_t = 5\cdot 10^{-3}$.

\begin{figure}
\begin{centering}
\includegraphics[width=.95\textwidth]{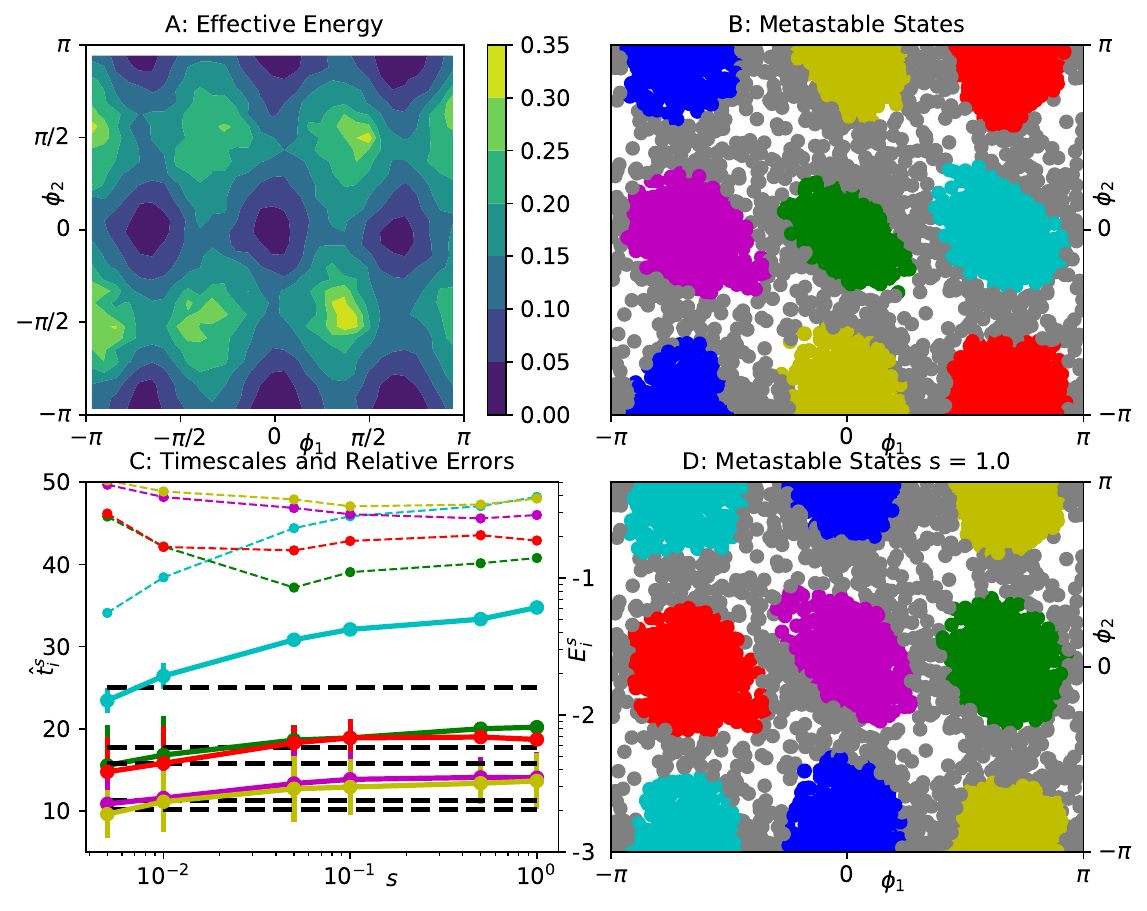}
\par\end{centering}
\centering{}\caption{Analysis of the effective dynamics of a prototypical five atom molecular system in the space of its dihedral angles $\phi_1,\, \phi_2$. A) Effective free energy in dihedral angle plane. B) Decomposition into six metastable states based on PCCA analysis of a gEDMD model with exact system parameters. Gray dots represent transition states where none of the memberships $\chi_j$ exceeds $0.6$. C) First five implied timescales $\hat{t}^s_i$ extracted from gEDMD models with KM estimators at various offsets $s$ (solid lines), compared to the gEDMD model with exact system parameters (dashed black lines). Errorbars were computed by bootstrapping. We also show the mean relative error $E^s_i$ given in \eqref{eq:relative_err_ev} (thin dashed lines, scale on the right, labeled by decadic logarithm). D) Decomposition into six metastable states based on PCCA analysis of a gEDMD model at offset $s = 1.0$. \label{fig:toy_molecule}}
\end{figure}

\noindent The system parameters are tuned in such a way that the two dihedral angles $\phi_1, \, \phi_2$, spanned by those five atoms, capture the slow dynamics of the system,  corresponding to transitions between six symmetrically arranged minima of the effective free energy in the $\phi_1$-$\phi_2$-plane, see Fig.~\ref{fig:toy_molecule} A. We run the reference gEDMD analysis using 100 two-dimensional Gaussian basis functions centered on a ten-by-ten grid, with~$\rho = 0.05$. Five dominant implied timescales between $t_1 \approx 25$ and $t_5 \approx 10$ (black lines in Fig.~\ref{fig:toy_molecule} C) are determined, and the six energy minima are recovered as metastable states by a PCCA analysis, as shown in Fig.~\ref{fig:toy_molecule} B.

\noindent We repeat the same experiment as for the previous example: using the same basis set of 100 Gaussians, we apply gEDMD with KM estimators \eqref{eq:gedmd_km}, for a series of offsets between $s = 5 \cdot 10^{-3}$ and $s = 1.0$. The corresponding timescale estimates $\hat{t}^s_i$ and the relative errors $E^s_i$ are also shown in Fig.~\ref{fig:toy_molecule}~C. The relative errors are generally larger than for the previous example, but still within an acceptable range of ten to forty percent. Applying the PCCA analysis to the model extracted at $s = 1.0$, we are still able to reproduce the correct metastable decomposition of reaction coordinate space. In summary, we can still conclude that the dominant spectrum is retained as long as the offset is at least an order of magnitude smaller than the fastest interesting timescale.

\subsection{Summary of Observations}
\label{ssec:sum_obs}
The previous two examples have indicated that, at least for reversible systems, and for a good selection of reaction coordinates, it is possible to use the KM estimators within the framework of gEDMD without disrupting the overall structure of the dominant spectrum. We state this observation as a conjecture, which will be investigated and made quantitative in future work:

\begin{Conjecture}
\label{conj:km_large_s}
Let $\xi$ be a good set of reaction coordinates for a reversible process $X_t$, which capture the slow part of the dynamics well (e.g. small projection errors $\delta_i$ as in Proposition \ref{prop:eigenvalue_error}). Let $s$ be an offset such that~$0 < s \ll t_K$. Then it is possible to recover the dominant spectrum by applying gEDMD with KM estimators, as in \eqref{eq:gedmd_km}, at least in a qualitative sense.
\end{Conjecture}

\noindent Note that Conjecture \ref{conj:km_large_s} is equivalent to saying that the effective dynamics with parameters \eqref{eq:Kramers_Moyal} will retain the dominant spectrum well in the setting described above. Also, the quality of a gEDMD approximation depends on the basis set. Conjecture \ref{conj:km_large_s} should be understood in the sense of using gEDMD with a powerful basis set, such that the basis set error does not play a major role.

\section{Underdamped Langevin Dynamics}
\label{sec:underdamped}
The third topic of this paper is to study the implications of the first two results for systems driven by underdamped Langevin dynamics. Let us recall that the theoretical analysis presented in sections \ref{sec:approx_spectrum} and \ref{sec:proofs} hinges on the reversible setting. Also, the numerical results on KM estimators shown in section \ref{sec:km_results} were using the reversible overdamped Langevin process \eqref{eq:overdamped_Langevin} (OL process from now on). However, a popular dynamical model, especially in molecular physics, are \emph{underdamped Langevin dynamics} (UL process) in position and momentum space $(q,\, p)$, where $q \in \Omega$ is the position of the system and $p \in \mathbb{R}^d$ is its momentum. The equations of motion are
\begin{eqnarray}
dq_{t} & = & p_{t}dt,\label{eq:Langevin_1}\\
dp_{t} & = & -\nabla V(q_{t})dt-\gamma p_{t}dt+\sqrt{2\gamma\beta^{-1}}dB_{t}.\label{eq:Langevin_2}
\end{eqnarray}
Here $V, \gamma, \beta$ have the same physical meaning of energy, friction, and inverse temperature, as in \eqref{eq:overdamped_Langevin} above. The invariant density of the process (\ref{eq:Langevin_1} - \ref{eq:Langevin_2}) factors as
\begin{equation}
\label{eq:mu_underdamped}
\mu(q, p) = \mu_P(p)\mu_Q(q) \propto \exp(-\frac{\beta}{2}p^T p)\exp(-\beta V(q)).
\end{equation}
 
\noindent The underdamped process is not reversible, so the error theory developed in this paper does not apply. However, it is well known that for large enough friction, if we observe the position coordinate $q_t$ only at every $s$-th step, for $s$ large, this time-rescaled process behaves like an OL process. This phenomenon is called the \emph{overdamped limit}, and gives rise to the following idea: if we select a reaction coordinate $\xi = \xi(q)$, which does not depend on momentum, then $\xi(q_{ks}),\, k=1, 2, \ldots$, observed at large offset $s$, behaves like a projected OL process. If we employ KM estimators or related approaches at the same offset within the context of gEDMD, we effectively compute models for the OL dynamics. Our error theory does apply to the reversible OL process, and if our Conjecture \ref{conj:km_large_s} is correct, we still get the dominant spectrum of these OL dynamics right by using KM formulae at a large offset. This idea will be illustrated by numerical examples in the following subsections.

\subsection{Projection and Re-scaling of the Underdamped Process}
We start by discussing the projection of the underdamped Langevin process by a map $\xi = \xi(q)$ depending only on the position coordinates. Oftentimes, one is not interested in quantities that explicitly depend on the momenta, which renders this a realistic setting (see \cite{DUONG2018} for an approach to model reduction which includes the momenta). Unfortunately, the coefficients of the effective dynamics \eqref{eq:def_projected_generator} are identically zero in this case, see also~\cite{Schutte:1999ab,Bittracher:2015aa}. To see this, it is readily checked using the definition of the generator that $\mathcal{L}\xi = p\cdot \nabla_q \xi$ and $\nabla \xi^T a \nabla \xi = 0$. Using the definition of the parameters \eqref{eq:projected_drift_diff}, the factorization \eqref{eq:mu_underdamped} of the invariant measure, and the fact that $\mu_P(p)$ is the density of a multivariate normal distribution with mean zero, we find that
\begin{align*}
b^\xi(z) &= \mathcal{P}(\mathcal{L}\xi)(z) \\
&= \int_{\Omega \times \mathbb{R}^d} \mathcal{L}\xi(q, p) \mu(q, p) \delta(\xi(q, p) - z) \,\mathrm{d}q \,\mathrm{d}p \\
&\propto \int_{\Omega \times \mathbb{R}^d} p\cdot \nabla_q \xi(q) \exp(-\beta V(q))\exp(-\frac{\beta}{2}p^T p) \delta(\xi(q) -z) \,\mathrm{d}q \,\mathrm{d}p \\
&= \left[\int_{\Omega} \exp(-\beta V(q)) \nabla_q^T \xi(q) \delta(\xi(q) - z) \,\mathrm{d}q \right] \cdot \left[ \int_{\mathbb{R}^d} p \exp(-\frac{\beta}{2}p^T p)\, \mathrm{d}p \right] \\
&= 0, \\
a^{\xi}(z) &= \mathcal{P}(\nabla \xi^T a \nabla \xi) = 0.
\end{align*}

\noindent How can one define a suitable effective dynamics in this case? As already mentioned above, if the friction $\gamma$ is sufficiently large, the positional component $q_t$ behaves like a reversible overdamped dynamics \eqref{eq:overdamped_Langevin} on long time scales. More precisely, if we observe the positions $q_{ks},\, k=1, 2, \ldots$ for an offset $s \gg \frac{1}{\gamma}$, then the statistics of this process will be approximately the same as those of an OL process $X_t$ in position space, observed at \emph{the same offset}. A particular pair of statistics to observe is given by the KM formulae \eqref{eq:Kramers_Moyal}. As a consequence, if the reaction coordinate $\xi$ captures the slow dynamics of the OL process $X_t$ well, and if in addition our Conjecture \ref{conj:km_large_s} is correct, then we can use the KM estimators on the underdamped data to build a good model (or a suitable effective dynamics) for the overdamped process $X_t$.

\noindent The following corollary provides a formal derivation of this argument, again in a qualitative sense, thus connecting the results of sections \ref{sec:approx_spectrum} and~\ref{sec:km_results}.

\begin{corollary} \label{cor:langevinApprox}
Let $X_t \in \Omega \subset \mathbb{R}^d$ denote the OL process \eqref{eq:overdamped_Langevin}. Moreover, assume $\xi$ is good reaction coordinate for $X_t$ (e.g. small projection errors $\delta_i$ in Proposition \ref{prop:eigenvalue_error}). Let $(q_t, p_t)$ denote the UL process on $\Omega \times \mathbb{R}^d$ with the same parameters as $X_t$.

(i) For $\gamma$ sufficiently large, the statistics of the positional component $q_t$ of the underdamped process, and of the overdamped process $X_t$, approximately agree at offsets much larger than $\frac{1}{\gamma}$.

(ii) If in addition Conjecture \ref{conj:km_large_s} is true, then application of the Kramers--Moyal estimator \eqref{eq:gedmd_km} at offsets $t_K > s > \frac{1}{\gamma}$ to $q_t$, allows to recover the dominant spectrum by means of gEDMD.
\end{corollary}

\begin{proof}
(i) Setting $\epsilon = \frac{1}{\gamma}$, and re-scaling the UL process (\ref{eq:Langevin_1} - \ref{eq:Langevin_2}) by $(q^c_t, p^c_t) = (q_{c\gamma t}, p_{c\gamma t}) = (q_{ct/\epsilon}, p_{ct/ \epsilon})$, the re-scaled equations of motion are
\begin{align*}
dq_{t}^c & = \frac{c}{\epsilon} p_{t}^c dt, \\
dp_{t}^c & = -\frac{c}{\epsilon}\nabla V(q_{t}^c)dt - \frac{c}{\epsilon^2} p_{t}^c dt+\frac{1}{\epsilon}\sqrt{2\beta^{-1}c} \, dB_{t}.
\end{align*}
For sufficiently large $\gamma$ and $c \geq 1$, Theorem 18.1 of \cite{Pavliotis:2008aa} implies that the law of $q_t^c$ is close to that of
\begin{equation*}
dQ^c_{t} = -c \nabla V(Q^c_{t}) dt + \sqrt{2\beta^{-1}c}\, dB_{t}.
\end{equation*}
We note that $Q^c_t$ is a time re-scaling of $X_t$ via $Q^c_t = X_{c\gamma t}$, hence we have for functions $f = f(q) \in L^2_{\mu_Q}(\Omega)$:
\begin{equation}
\label{eq:expectations_langevin_long_timescales}
\mathbb{E}^{\mu_Q}\left[f(q_{c\gamma t})\right] \approx \mathbb{E}^{\mu_Q}\left[f(X_{c\gamma t})\right].
\end{equation}

(ii) In particular, the approximate equality \eqref{eq:expectations_langevin_long_timescales} applies directly to estimator \eqref{eq:gedmd_km} for $s > \frac{1}{\gamma}$. If Conjecture \ref{conj:km_large_s} is true, the resulting gEDMD models will qualitatively retain the leading timescales of the process $X_t$.
\end{proof}

\subsection{Langevin Toy Model}
\label{ssec:Langevin_toy_model}
\begin{figure}
\begin{centering}
\includegraphics[width=.7\textwidth]{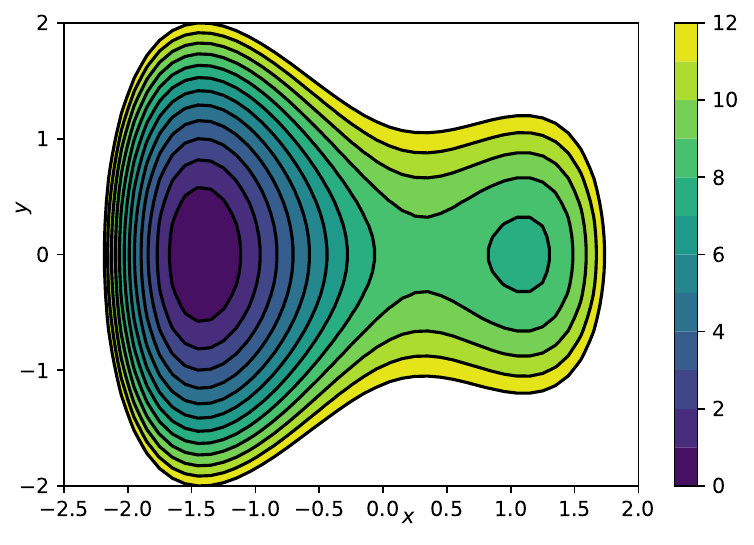}
\par\end{centering}
\centering{}\caption{Two-dimensional model potential \eqref{eq:2d_toy_potential}. \label{fig:2d_toy_potential}}
\end{figure}

\begin{figure}
\begin{centering}
\includegraphics[width=.95\textwidth]{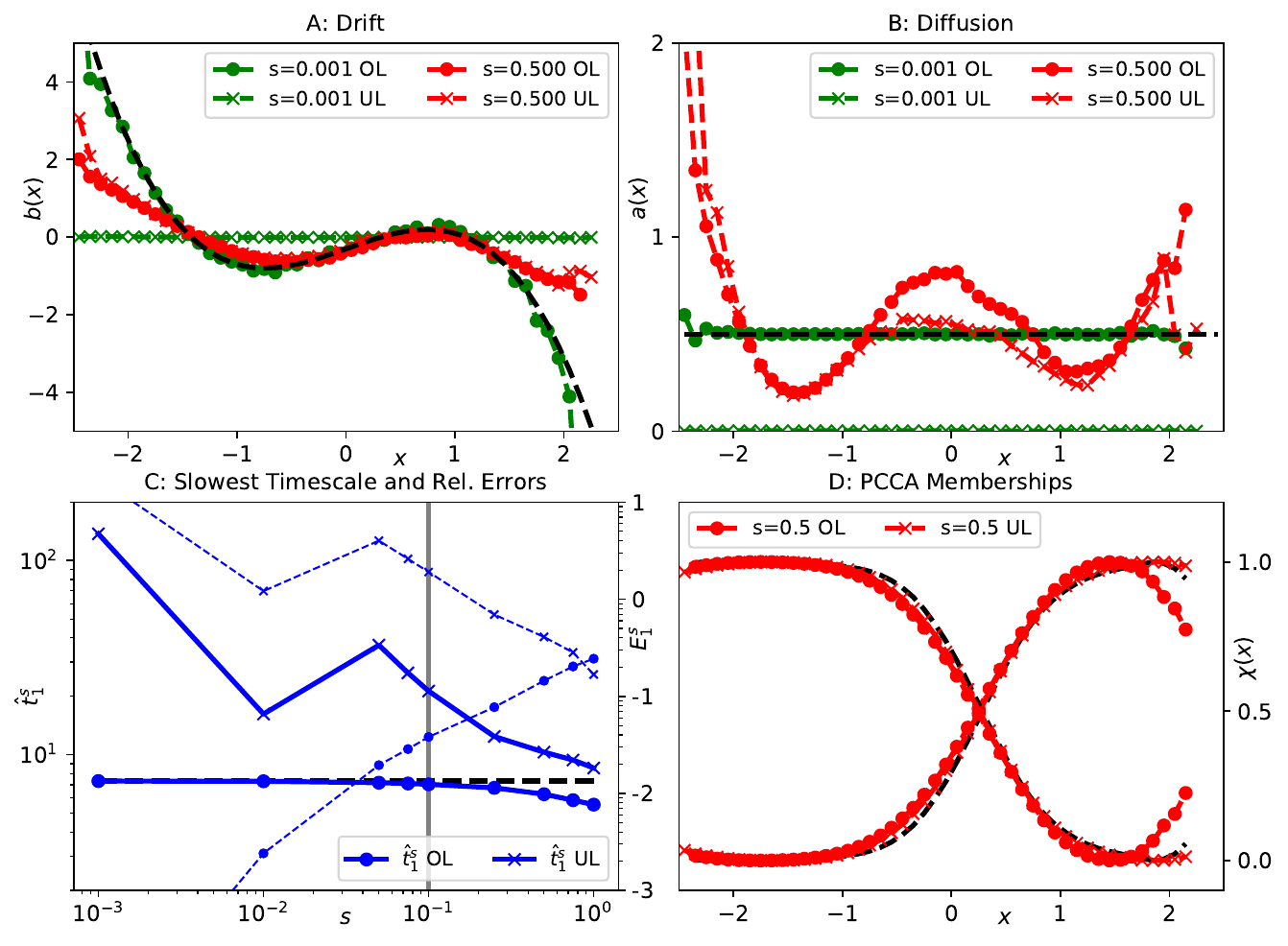}
\par\end{centering}
\centering{}\caption{Analysis of effective dynamics on the $x$-coordinate of the two-dimensional
toy potential Eq. (\ref{eq:2d_toy_potential}). A: Numerical estimates of effective drift for different values of the offset $s$ using both OL data (dots) and UL data (crosses). B: The same for the effective diffusion. C: Leading implied timescale $\hat{t}_1^s$ obtained from gEDMD models of the projected OL data (dots) and UL data (crosses), as a function of $s$, compared to the reference value in black. The reference was extracted from a gEDMD model using exact system parameters. We also show the relative error $E^s_i$ for both data sets (thin dashed lines, scale on the right, labeled by decadic logarithm). The vertical gray line indicates the critical relaxation time $\frac{1}{\gamma}$. D: PCCA memberships extracted from gEDMD models at offset $s = 0.5$ for both data sets, compared to the reference gEDMD model in black. \label{fig:langevin_2d}}
\end{figure}

To illustrate the ideas outlined in the previous section, we first consider another two-dimensional toy potential
\begin{equation}
\label{eq:2d_toy_potential}
V(x, y) = 3x^4 - 5x^2 + 1.5x + 3y^2,
\end{equation}
shown in Fig. \ref{fig:2d_toy_potential}. For $\gamma = 10, \, \beta = 0.4$, we generate data of both the UL and OL processes for this potential, each data set comprising $M = 10^7$ points at integration time step $\Delta_t = 10^{-3}$. For both dynamical models, the slowest transition in this energy landscape is the crossing of the barrier around $x = 0$, thus $\xi(x, y) = x$ is a suitable reaction coordinate in both instances. The associated implied timescale in the OL model is $t_1 \approx 7.3$. The effective drift in the overdamped case simply corresponds to the $x$-derivative of the $x$-dependent part the potential, while the effective diffusion remains constant, see the black dashed lines in Figs. \ref{fig:langevin_2d} A and B.

\noindent We first illustrate the relationship between the overdamped limit and estimates of drift and diffusion by means of KM formulae. In Figures \ref{fig:langevin_2d} A and B, we show histogrammed estimates of the KM expressions \eqref{eq:Kramers_Moyal}, using both the OL and UL simulation data. We see that for a small offset $s = 10^{-3}$, UL estimates are almost zero, while the OL data lead to estimates close to the analytical values. Both findings are as expected. For a much larger offset $s = 0.5$, however, both estimates are significantly different from what we find for small offsets, and most importantly, the overdamped and underdamped estimates agree well. This in line with our qualitative argument in Corollary \ref{cor:langevinApprox} (i). 

\noindent We apply gEDMD to both data sets, using the KM estimators as in \eqref{eq:gedmd_km}, for a series of offsets between $s = 10^{-3}$ and $s = 1.0$. The basis set is comprised of fifteen Gaussian basis functions \eqref{eq:gaussian_fct}, centered uniformly between $x = -2.5$ and $x = 2.5$, each of width $\rho = 0.1$. From all of these models, we extract the slowest implied timescale $\hat{t}_1^s$ and the relative errors $E^s_1$, and present the results in Fig. \ref{fig:langevin_2d} C. As expected, OL estimates are highly accurate at small offsets. For large offsets, these estimates remain within the same error margin that we have seen before. Estimates based on the UL data, however, are far off for small offsets (as the coefficients of the corresponding dynamics are essentially zero), but once the offset exceeds the critical value $s > \frac{1}{\gamma} = 0.1$, they are about as accurate as those based on OL data. We also verify in Fig.~\ref{fig:langevin_2d} D that for both data sets, the two metastable states of the effective dynamics can be correctly identified by applying PCCA to the eigenvectors of the gEDMD model at a large offset $s = 0.5$. 

\noindent This example shows that upon increasing the offset it is possible to find a sweet spot where $s$ is larger than the critical relaxation time $\frac{1}{\gamma}$, but smaller than the slowest timescale $t_1$, such that meaningful effective dynamics for the UL process along $x$ can be defined in this regime.

\subsection{Alanine Dipeptide}
\label{ssec:alanine_dipeptide}
\begin{figure}
\begin{centering}
\includegraphics[width=.95\textwidth]{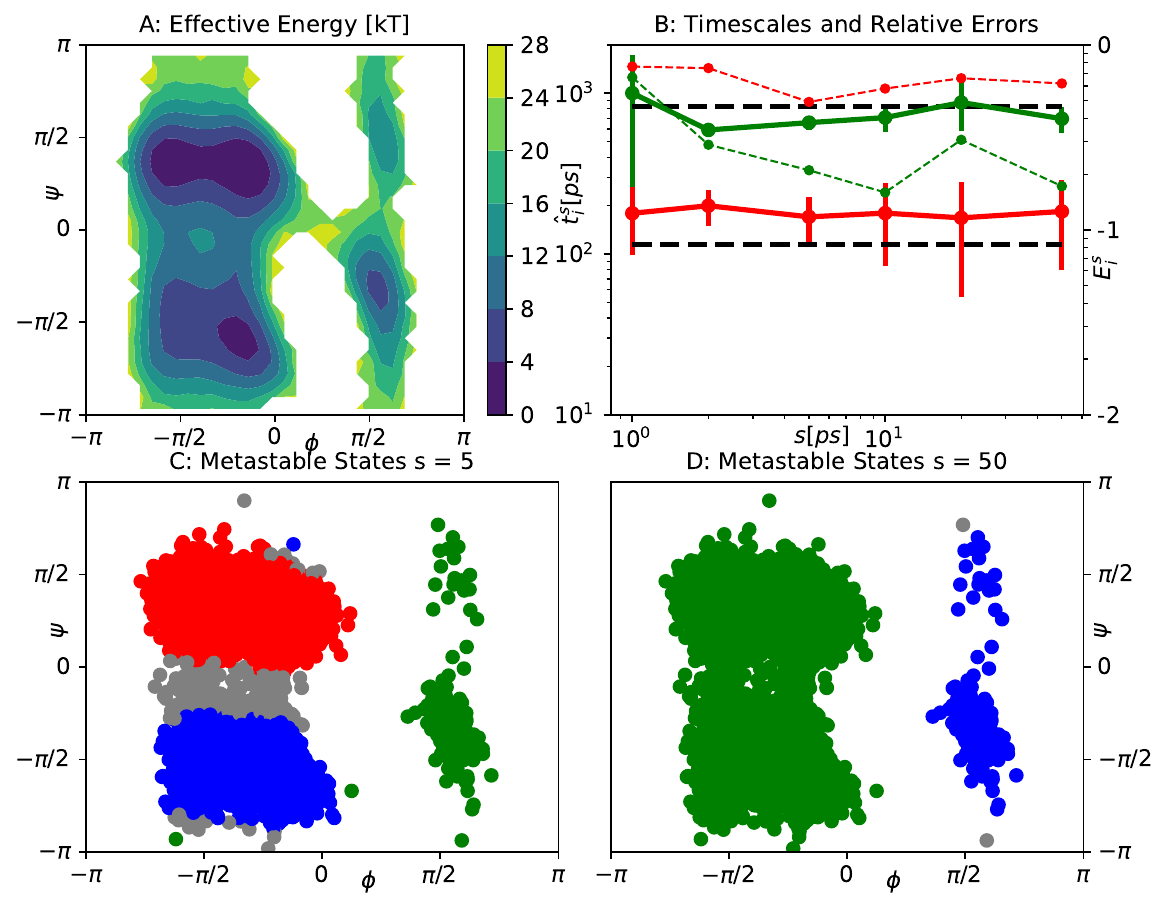}
\par\end{centering}
\centering{}\caption{Analysis of effective dynamics for alanine dipeptide in the space of its backbone dihedral angles $\phi,\, \psi$.
A: Effective free energy of the original simulation data in the $\phi-\psi$-plane. Metastable states correspond to the two deep minima on the left, and the shallow minimum on the right. B: Slowest two timescales $\hat{t}_{1}^s, \, \hat{t}_2^s$ computed by gEDMD models at various offsets $s$, compared to the reference values in black. Errorbars were estimated by bootstrapping. We also show the mean relative error $E_i^s$ \eqref{eq:relative_err_ev} (thin dashed lines, scale on the right, labeled by decadic logarithm). C: Metastable decomposition into three states determined by applying PCCA to the eigenfunctions of the gEDMD model at $s = 5\,\mathrm{ps}$. Gray dots represent transition states where none of the memberships $\chi_j$ exceeds 0.6. D: The same for $s = 50\,\mathrm{ps}$, but using only two states.  \label{fig:ala2}}
\end{figure}

Our final numerical example is a more complex data set, namely molecular dynamics simulations of the alanine dipeptide. The data we use is the same as in Refs.~\cite{Nuske2017,WANG2019}, consisting of $M = 10^6$ points at $1 \mathrm{ps}$ time spacing. The dynamical model is the UL process as in section \ref{ssec:Langevin_toy_model}, using the AMBER 99 molecular force field \cite{LINDORFF2010} as potential $V$. Friction is set to $\gamma = 0.1 \mathrm{ps}^{-1}$, and $\beta$ is derived from the temperature $T = 300 \mathrm{K}$ via $\beta^{-1} = k_B T$, using the Boltzmann constant $k_B$. As is well-known from numerous previous studies, the slow dynamics of alanine dipeptide can be represented well in the space of backbone dihedral angles $\phi,\,\psi$.  Figure \ref{fig:ala2} A shows the effective free energy landscape in the space of these reaction coordinates. Three major minima can be identified, which correspond to three metastable states of the full dynamics. The corresponding transition timescales $t_1 \approx 1 \, \mathrm{ns}$ and $t_2 \approx 0.1 \, \mathrm{ns}$ are indicated by the black lines in Fig.~\ref{fig:ala2} B.

\noindent For gEDMD, we employ a similar basis set as in section \ref{ssec:toy_molecule}, comprised of 100 periodic Gaussians, centered on a ten-by-ten grid between $-2.5$ and $2.5$, with bandwidth $\rho = 0.05$. We compute gEDMD models for a range of offsets between $s = 1 \mathrm{ps}$ and $s = 50 \mathrm{ps}$. The resulting first two implied timescales $\hat{t}_1^s,\, \hat{t}_2^s$ and the corresponding relative errors $E_i^s$ are shown in Fig. \ref{fig:ala2} B. These results confirm the findings of the previous examples, as both leading timescales are roughly reproduced for all offsets considered. The relative error $E_2^s$ for the second timescale is generally larger than what we observed in the previous examples, but it is still acceptable in the context of this example. An interesting observation is that for small offsets, we are able to fully recover both slow processes by a PCCA analysis. The corresponding state decomposition is shown in Fig.~\ref{fig:ala2} C. For a larger offset which is comparable to the timescale $t_2$, the faster transition within the left part of the plane appears to be blurred out. However, applying PCCA with only two metastable states still recovers the slowest process, as indicated by the decomposition in Fig. \ref{fig:ala2} D.

\section{Precise statements on spectral properties and their proofs}
\label{sec:proofs}
In this section, we provide detailed proofs of the spectral approximation results outlined in section \ref{sec:approx_spectrum}.

\subsection{Form Domain}
\label{ssec:form_domain}
\noindent We consider an open domain $\Omega\subset \mathbb{R}^d$, and make the assumption of uniform ellipticity \eqref{eq:uniform_ellipticity}. The generator $\mathcal{L}$ in \eqref{eq:generator_differential_op} can be defined initially on the set of smooth functions or smooth and compactly supported functions. The form domain $\mathbb{V}_\mu$ can then be obtained as the closure of this initial domain with respect to the \emph{Dirichlet norm}
\begin{equation}
\label{eq:dirichlet_norm_2}
\|\psi \|^2_1 = \| \psi \|^2_{L^2_\mu} + Q(\psi, \psi).
\end{equation}
We note that on domains with a boundary, the choice of initial domain has an impact on the boundary conditions.
%, but we do not go into details of this matter here.
In addition, we can restrict all function spaces to the orthogonal complement of the constant one function without explicitly changing the notation. 

\noindent The assumption of uniform ellipticity implies that $\mathbb{V}_\mu$ is also a Hilbert space if equipped with the energy norm
\begin{equation}
\label{eq:energy_norm_2}
\|\cdot \|_\mathcal{Q} = \mathcal{Q}(\cdot, \cdot)^{1/2}.
\end{equation}
In many cases of practical interest, the energy and Dirichlet norms are equivalent \cite{BAKRY2013}. We state this as an assumption and generally use the energy norm on $\mathbb{V}_\mu$ in what follows:

\noindent \textbf{Assumption 1:} The Dirichlet norm \eqref{eq:dirichlet_norm_2} and energy norm \eqref{eq:energy_norm_2} are equivalent.
.

\subsection{Solution Operator and Discrete Spectrum}
The solution operator $\mathcal{T}:\,\mathbb{V}_\mu \rightarrow \mathbb{V}_\mu$ associated to $\mathcal{Q}$ is defined by
\begin{equation}
\label{eq:definition_solution_op}
\mathcal{Q}(\mathcal{T}\psi, \tilde{\psi}) = \innerprod{\psi}{\tilde{\psi}}_\mu \quad \forall \tilde{\psi} \in \mathbb{V}_\mu.
\end{equation}

\noindent \textbf{Assumption 2}: The solution operator is compact on $\mathbb{V}_\mu$ with norm~$\| \cdot \|_\mathcal{Q}$.

\noindent As a consequence of this assumption, the generator possesses a complete set of eigenfunctions $\psi_i \in \mathbb{V}_\mu$ with eigenvalues $\kappa_i, \, i = 1, 2, \ldots$, which are given as reciprocals of the eigenvalues of $\mathcal{T}$. There is a number of well-known settings where Assumption 2 can be shown to hold. These include bounded Lipschitz domains with Dirichlet, Neumann or periodic boundary conditions, as well as overdamped Langevin dynamics with a potential satisfying suitable growth conditions on the potential $V$~\cite{BAKRY2013}.

\subsection{Coarse-Grained Generator and its Spectrum}
Analogously to section \ref{ssec:form_domain}, the projected generator $\mathcal{L}^\xi$ can be defined initially on the set of smooth or smooth and compactly supported functions in $L^2_\nu$, which is an infinite-dimensional subspace of $L^2_\mu$ \cite{Zhang:2016aa}. Using the effective quadratic form $\mathcal{Q}^\xi$ \eqref{eq:effective_quad_form}, the effective form domain $\mathbb{V}_\nu$ can be defined again by completion of the initial domain with respect to the corresponding Dirichlet norm $\|\phi\|^2_{1, \xi} = \|\phi\|^2_{L^2_\nu} + \mathcal{Q}^\xi(\phi, \phi)$. Due to the relations \cite{Zhang:2016aa}
\begin{align}
\label{eq:exchange_inner_products}
\innerprod{\phi}{\tilde{\phi}}_\nu &= \innerprod{\phi\circ \xi }{\tilde{\phi}\circ \xi}_\mu, & \mathcal{Q}^\xi(\phi, \tilde{\phi}) &= \mathcal{Q}(\phi\circ \xi, \tilde{\phi}\circ \xi),
\end{align}
Assumption 1 carries over to the effective Dirichlet norm and the effective energy norm. In order to proceed from this point, some care needs to be taken with regard to the coarse graining map:

\noindent \textbf{Assumption 3:} The coarse graining map $\xi$ is such that the effective form domain $\mathbb{V}_\nu$ is a subspace of~$\mathbb{V}_\mu$.

\begin{remark} Assumption 3 will hold in many cases of practical interest, which include the projection of a periodic domain onto a lower-dimensional periodic domain by a function $\xi$ which respects the periodic boundary conditions. As a negative example, however, consider an SDE on a rectangle in two dimensional space, with absorbing boundary conditions. Choose the coarse graining function as the projection onto the first coordinate axis. The form domain and the effective form domain will be given as first order Sobolev spaces with zero boundary conditions, but the effective form domain is not contained in the full form domain as its elements do not vanish on parts of the full boundary which are parallel to the first coordinate axis. 
\end{remark}

\noindent If Assumption 3 holds, it makes sense to define the $Q$-orthogonal projection from $\mathbb{V}_\mu$ onto $\mathbb{V}_\nu$, which we denote by $\mathcal{P}_\mathcal{Q}$. As a first main result, we show that assumptions 1, 2 and 3 are sufficient to ensure that the spectrum of $\mathcal{L}^\xi$ is also discrete:

\begin{proposition}
\label{prop:discrete_spectrum_lxi}
If assumptions 1, 2, and 3 hold, the spectrum of $\mathcal{L}^\xi$ is discrete.
\end{proposition}
\begin{proof}
By assumption 1, the effective solution operator $\mathcal{T}^\xi$ is uniquely defined on $\mathbb{V}_\nu$ by
\begin{equation*}
\mathcal{Q}^\xi(\mathcal{T}^\xi \phi, \tilde{\phi}) = \innerprod{\phi}{\tilde{\phi}}_\nu \quad \forall \tilde{\phi} \in \mathbb{V}_\nu.
\end{equation*}
It is readily checked that $\mathcal{T}^\xi = \mathcal{P}_\mathcal{Q}\mathcal{T}\mathcal{P}_Q$. Since $\mathcal{P}_\mathcal{Q}$ is bounded and $\mathcal{T}$ is compact by assumption 2, so is $\mathcal{T}^\xi$. Hence its spectrum (and correspondingly that of $\mathcal{L}^\xi$) is discrete.
\end{proof}

\subsection{Approximation Result}

Next, we derive a bound on the eigenvalue error of $\mathcal{L}^\xi$ in terms of the energy-norm approximation error of the dominant eigenfunctions of $\mathcal{L}$. The idea is to apply classical Galerkin error estimates to a sequence of finite-dimensional subspaces $\mathbb{W}^h$ in $\mathbb{V}_\nu$, and exploit that these provide approximations to the (nonzero) eigenvalues of both $\mathcal{L}$ and $\mathcal{L}^\xi$.

\begin{proposition}
\label{prop:eigenvalue_error}Denote the projection error of eigenfunction $\psi_{i}$ with respect to the energy-norm by $\delta_i = \|\mathcal{P}^\perp_\mathcal{Q} \psi_i\|_\mathcal{Q}$. The relative error between the $i$-th eigenvalues of $\mathcal{L}^{\xi}$ and $\mathcal{L}$ is bounded by
\begin{equation}
\frac{\omega_{i}-\kappa_{i}}{\omega_{i}} \leq \left[1 + \max_{j=1, \ldots, i-1} \frac{\omega_j^2 \kappa_{i}^2}{(\omega_j - \kappa_i)^2}\| (\mathcal{I} - \mathcal{P}_\mathcal{Q})\mathcal{T}\|_\mathcal{Q}^2 \right] \delta_i^2. \label{eq:error_bound_rates}
\end{equation}
\end{proposition}

\begin{proof}
We consider a sequence $\mathbb{W}^h$ of finite-dimensional subspaces in the reduced space $\mathbb{V}_\nu$. For every $\mathbb{W}^h$, the $\mathcal{Q}$-orthogonal projection from $\mathbb{V}_\mu$ onto $\mathbb{W}^h$ is labeled by $\mathcal{P}^h_\mathcal{Q}$, and the corresponding $\mathcal{Q}^\xi$-orthogonal projection in $\mathbb{V}_\nu$ is called $\mathcal{P}^h_{\mathcal{Q},\xi}$. We assume that all $\mathbb{W}^h$ contain the projections $\mathcal{P}_\mathcal{Q} \psi_j$ for $j = 1\, \ldots, i$, and satisfy the following approximability condition in $\mathbb{V}_\nu$ (which holds in any separable Hilbert space):
\begin{equation}
\label{eq:pointwise_convergence_ph}
\lim_{h\rightarrow 0} \|(\mathcal{I} - \mathcal{P}^h_{\mathcal{Q}, \xi}) \phi \|_{\mathcal{Q}^\xi} = 0 \quad \forall \phi \in \mathbb{V}_\nu.
\end{equation}
Note that the spaces $\mathbb{W}^h$ are just of auxiliary nature; their sole purpose being to reconcile the theory of~\cite{KYAZEV2006}---that uses finite-dimensional approximation spaces---with our infinite-dimensional~$\mathbb{V}_{\nu}$. There is no need to specify them in detail, and $h$ only serves as a formal parameter here. Using \eqref{eq:min_max_L_Lxi}, and by applying \cite[Theorem 3.2]{KYAZEV2006} to the approximation of $\mathcal{L}$, we find that for any of the approximation spaces $\mathbb{W}^h$:
\begin{align*}
\frac{\omega_{i}-\kappa_{i}}{\omega_{i}} \leq \frac{\hat{\omega}_{i}-\kappa_{i}}{\hat{\omega}_{i}} &\leq \left[1 + \max_{j=1, \ldots, i-1} \frac{\hat{\omega}_{j}^2 \kappa_{i}^2}{(\hat{\omega}_{j} - \kappa_i)^2} \| (\mathcal{I} - \mathcal{P}^h_\mathcal{Q})\mathcal{T}\|_Q^2 \right] \|(\mathcal{I} - \mathcal{P}^h_\mathcal{Q}) \psi_i \|_Q^2 \\
&= \left[1 + \max_{j=1, \ldots, i-1} \frac{\hat{\omega}_{j}^2 \kappa_{i}^2}{(\hat{\omega}_{j} - \kappa_i)^2} \| (\mathcal{I} - \mathcal{P}^h_\mathcal{Q})\mathcal{T}\|_Q^2 \right] \delta_i^2,
\end{align*}
where $\hat{\omega}_i$ are the Ritz values associated to $\mathbb{W}^h$, see section \ref{sec:MSM_VAC}. The last equality holds because $\mathcal{P}_Q \psi_i \in \mathbb{W}^h$ by assumption. It remains to study the pre-factor in the limit of $h \rightarrow 0$. We conclude from \eqref{eq:pointwise_convergence_ph} that $\hat{\omega}_{i} \rightarrow \omega_i$ as $h \rightarrow 0$, which already yields the first term in the pre-factor. Regarding the second term, we first observe that $\mathcal{P}^h_\mathcal{Q} = \mathcal{P}^h_{\mathcal{Q}, \xi} \mathcal{P}_\mathcal{Q}$, since for any $\psi \in \mathbb{V}_\mu$ and $\phi \in \mathbb{W}^h$:
\begin{align*}
\innerprod{\psi - \mathcal{P}^h_{\mathcal{Q}, \xi} \mathcal{P}_\mathcal{Q} \psi }{\phi}_\mathcal{Q} 
&= \innerprod{(\mathcal{P}_\mathcal{Q} + \mathcal{P}_\mathcal{Q}^\perp) \psi - \mathcal{P}^h_{\mathcal{Q}, \xi} \mathcal{P}_\mathcal{Q} \psi }{\phi}_\mathcal{Q} \\
&= \innerprod{\mathcal{P}_\mathcal{Q} \psi - \mathcal{P}^h_{\mathcal{Q}, \xi} \mathcal{P}_\mathcal{Q} \psi }{\phi}_\mathcal{Q}
= \innerprod{(\mathcal{I} - \mathcal{P}^h_{\mathcal{Q}, \xi}) \mathcal{P}_\mathcal{Q} \psi}{\phi}_{\mathcal{Q}^\xi} = 0.
\end{align*}
The third equality is due to \eqref{eq:exchange_inner_products}, while the last equality is due to the definition of $\mathcal{P}^h_{\mathcal{Q}, \xi}$. With this result, and using pointwise convergence \eqref{eq:pointwise_convergence_ph} of $\mathcal{P}^h_{\mathcal{Q},\xi}$ to the identity in $\mathbb{V}_\nu$, we conclude for any $\psi \in \mathbb{V}_\mu$:
\begin{align*}
\lim_{h\rightarrow 0} \|(\mathcal{I} - \mathcal{P}^h_\mathcal{Q})\psi - \mathcal{P}^\perp_\mathcal{Q} \psi \|_\mathcal{Q}  &= \lim_{h\rightarrow 0} \|(\mathcal{P}_\mathcal{Q} - \mathcal{P}^h_\mathcal{Q}) \psi \|_{\mathcal{Q}} \\
&= \lim_{h\rightarrow 0} \|(\mathcal{I} - \mathcal{P}^h_{\mathcal{Q},\xi}) \mathcal{P}_Q\psi \|_{\mathcal{Q}^\xi} = 0.
\end{align*}
Combing this pointwise convergence of $\mathcal{I} - \mathcal{P}^h_\mathcal{Q}$ towards $\mathcal{P}_Q^\perp$ in $\mathbb{V}_\mu$, and compactness of the solution operator, we have
\begin{equation*}
\lim_{h \rightarrow 0} (\mathcal{I} - \mathcal{P}^h_\mathcal{Q})\mathcal{T} = (\mathcal{I} - \mathcal{P}_\mathcal{Q})\mathcal{T}
\end{equation*}
with respect to the operator norm. This establishes the second term in the pre-factor and hence the proposition.
\end{proof}

\noindent To make the bound from the previous result more accessible, we will bound the term involving the operator $\mathcal{T}$.
\begin{corollary} \label{cor:simpler_bound}
In the setting of Proposition~\ref{prop:eigenvalue_error}, the following error bound holds:
\begin{equation}
\frac{\omega_{i}-\kappa_{i}}{\omega_{i}} \leq \left[1 + \frac{ \kappa_{i}^2}{ \kappa_1^2} \max_{j=1, \ldots, i-1} \frac{\omega_j^2}{(\omega_j - \kappa_i)^2} \right] \delta_i^2. \label{eq:error_bound_rates2}
\end{equation}
\end{corollary}
\begin{proof}
As already mentioned, it follows from \eqref{eq:definition_solution_op} that
\[
\mathcal{T} \psi_i = \frac{1}{\kappa_i} \psi_i \qquad \forall i \geq 1.
\]
Since the $\psi_i$ are $\mathcal{Q}$-orthogonal, it follows that $\| \mathcal{T} \|_{\mathcal{Q}} = \frac{1}{\kappa_1}$. Further, since $\mathcal{P}_Q$ is a $\mathcal{Q}$-orthogonal projection, so is $\mathcal{I} - \mathcal{P}_Q$, and thus both have $\mathcal{Q}$-norm at most one. It follows that $\| (\mathcal{I} - \mathcal{P}_Q)\mathcal{T} \|_{\mathcal{Q}} \le \frac{1}{\kappa_1}$, and hence the claim.
\end{proof}

\begin{remark}
\label{rem:QnormHnorm}
Due to the equivalence of the norms $\|\cdot \|_\mathcal{Q}$ and $\|\cdot \|_{1}$, the bound in Proposition \ref{prop:eigenvalue_error} also applies to the latter norm. Let $\mathcal{P}_1$ denote the orthogonal projection onto $\mathbb{V}_\nu$ with respect to the Dirichlet norm. Then
\begin{equation*}
\|\mathcal{P}_\mathcal{Q}^\perp \psi_i \|_\mathcal{Q} \leq \|\mathcal{P}_1^\perp \psi_i \|_\mathcal{Q} \leq C \|\mathcal{P}_1^\perp \psi_i \|_{1}.
\end{equation*}
\end{remark}

\subsection{Comments}

The bound from Corollary~\ref{cor:simpler_bound} allows some interpretation. Note that in \eqref{eq:error_bound_rates2} the ratio of the eigenvalues $\kappa_1,\kappa_i$ of $\mathcal{L}$ and the relative difference of eigenvalues $\kappa_i$ and $\omega_j$ of $\mathcal{L}$ and $\mathcal{L}^{\xi}$ (respectively) play a role. Let us fix $i$ and assume that the first $i$ timescales are comparable, meaning that $\kappa_i/\kappa_1$ is a moderate number. Let us further assume that the squared projection errors $\delta_{\ell}$ are all much smaller than the \emph{relative eigenvalue gaps}, i.e.,
\[
\delta_{\ell}^2 \ll \min\bigg \{1,\frac{\kappa_{j} - \kappa_{j-1}}{\kappa_j} \bigg \} \qquad \forall j,\ell = 1,\ldots,i.
\]
Then, inductively from $j=1$ to $j=i$ it follows that $(\omega_j - \kappa_j) / \omega_j \approx \mathcal{O}(\delta_j^2)$, since we find that $\omega_j \approx \kappa_j$, thus the relative differences between $\kappa_{\ell}$ and $\omega_j$ are large for $j\neq \ell$, which makes the term in the square brackets on the right-hand side of \eqref{eq:error_bound_rates2} moderate. Hence, in this situation the relative error of the $i$-th timescale's approximation is effectively governed by the projection error~$\delta_i^2$.

\noindent We also note that the bound assumed in~\cite[Theorem~2]{Zhang2017} implies our bound up to another multiplicative constant:

\begin{lemma}
\label{lem:statement_adn_lemma}
Let us consider the diffusion~(\ref{eq:ito_sde}), satisfying~(\ref{eq:uniform_ellipticity}), on a bounded domain with smooth boundary and reflecting boundary conditions.
If $\|\mathcal{L}\mathcal{P}^{\perp} \psi_{i}\|_{L_{\mu}^{2}}\leq\delta_{i}$, then there is a~$C>0$ such that the assumptions of Proposition~\ref{prop:eigenvalue_error} are satisfied with
\begin{eqnarray}
\|\mathcal{P}^{\perp}_\mathcal{Q} \psi_{i}\|_{\mathcal{Q}} & \leq & C \delta_{i}.\label{eq:statement_adn_lemma}
\end{eqnarray}
\end{lemma}

\begin{proof}
We use the weighted Sobolev spaces $H^k_\mu$ \cite{CeDu18}. Since the spatial domain is bounded and~$\mu$ is smooth, it is a weight in Muckenhoupt class, cf.~\cite[Eq.~(1.2)]{CeDu18}. The regularity conditions \cite[Eqs.~(2.2)--(2.4)]{CeDu18} are satisfied by assumption.
Now, the result follows directly from the weighted Agmon--Douglis--Nierenberg estimate~\cite[Theorem 2.4]{CeDu18} giving~$\|u\|_{H^2_{\mu}} \le C\| \mathcal{L}u \|_{L^2_{\mu}}$ for some~$C>0$ independent of~$u$. Using assumption 1 and \eqref{eq:uniform_ellipticity}, we have

\begin{align*}
\|\mathcal{P}^{\perp}_\mathcal{Q} \psi_{i}\|_{\mathcal{Q}} &\leq \|\mathcal{P}^{\perp} \psi_{i}\|_{\mathcal{Q}} \leq C\|\mathcal{P}^{\perp} \psi_{i}\|_{H_{\mu}^{1}}  \leq C \|\mathcal{P}^{\perp} \psi_{i}\|_{H_{\mu}^{2}} \\
 & \leq C \|\mathcal{L}\mathcal{P}^{\perp} \psi_{i}\|_{L_{\mu}^{2}} \leq C \delta_i.
\end{align*}
concluding the proof.
\end{proof}

\section{Conclusions}
\label{sec:conclusion}
We have investigated the approximation of high-dimensional diffusion processes by effective dynamics defined on the lower-dimensional space of reduced variables. For reversible diffusions, a new relative error bound for the approximation of eigenvalues of the infinitesimal generator by the eigenvalues of a reduced generator was proved. Our bound shows that a small projection error of the corresponding generator eigenfunctions, measured by the energy norm, is sufficient for small eigenvalue errors.

\noindent In addition, we have presented numerical examples regarding the data-driven estimation of the eigenvalues of projected generators by means of the gEDMD method. If the full system parameters are unknown, they need to be approximated, for example using Kramers--Moyal formulae. We have presented numerical examples that for reversible systems, and good reaction coordinates, the resulting spectral estimates seem to remain stable across a long range of time windows in the KM-estimators.

\noindent Finally, we have suggested a strategy to define meaningful effective equations for underdamped Langevin dynamics on a subset of its position space. Exploiting the overdamped limit, and using KM estimators at large time windows, we can effectively model a projected overdamped process on the same domain. Numerical examples have confirmed the feasibility of this approach.

Future work will focus on providing a theoretical foundation for the observations stated in Conjecture~\ref{conj:km_large_s}. Also, the relation between the positional coordinate of an (underdamped) Langevin process on long timescales, and the corresponding  overdamped Langevin equation (cf. Corollary~\ref{cor:langevinApprox} (i)), needs to be analyzed in more detail.

%%%%%%%%%%%%%%%%%%%%%%%%%%%%%%%%%%%%%%%%%%
\section*{Funding and Acknowledgments} 
This research was funded by Rice University Academy of Fellows (F.N.), by Deutsche Forschungsgemeinschaft (DFG) through CRC 1114 "Scaling Cascades in Complex Systems", Project Number 235221301, project A01 (P.K.), by the Einstein Foundation Berlin (C.C.), by the National Science Foundation (Nos. CHE-1265929, CHE- 1738990, and PHY-1427654) and the Welch Foundation (No. C-1570) (F.N., L.B., C.C.).

\noindent We would like to thank Ralf Banisch, Stefan Klus, Thomas Swinburne, Tony Leli\`{e}vre and Frank No\'{e}. We are also indebted to the Institute for Pure and Applied Mathematics (IPAM) for hosting the long program "Complex High-Dimensional Energy Landscapes".

\appendix

\section{Lemon Slice Potential: Effective Drift and Diffusion Coefficients}

\label{sec:appendix_lemon_slice}

Here, we calculate the effective drift and diffusion if the lemon slice potential Eq.~(\ref{eq:lemon_slice}) is projected onto the polar angle $\xi(x,y)=\varphi(x,y)$. We start by expressing the generator $\mathcal{L}$ in polar coordinates. The Laplacian operator in polar coordinates is
\begin{equation*}
\Delta f = \frac{\partial^{2}f}{\partial r^{2}}+\frac{1}{r^{2}}\frac{\partial^{2}f}{\partial\varphi^{2}}+\frac{1}{r}\frac{\partial f}{\partial r}.
\end{equation*}
Moreover, it follows from the chain rule and the definition of $\xi$ that for any function $f$,
\begin{align*}
\frac{\partial f}{\partial x} & = \cos\varphi\frac{\partial f}{\partial r}-\frac{\sin\varphi}{r}\frac{\partial f}{\partial\varphi},\\
\frac{\partial f}{\partial y} & =  \sin\varphi\frac{\partial f}{\partial r}+\frac{\cos\varphi}{r}\frac{\partial f}{\partial\varphi}.
\end{align*}
The generator in polar coordinates then becomes
\begin{equation*}
\mathcal{L} =  -\left[\frac{\partial V}{\partial r}\frac{\partial}{\partial r}+\frac{1}{r^{2}}\frac{\partial V}{\partial\varphi}\frac{\partial}{\partial\varphi}\right]+\frac{1}{\beta}\left[\frac{\partial^{2}}{\partial r^{2}}+\frac{1}{r^{2}}\frac{\partial^{2}}{\partial\varphi^{2}}+\frac{1}{r}\frac{\partial}{\partial r}\right].
\end{equation*}
Applying the generator to the reaction coordinate $\xi$, only one of the terms above is non-zero, resulting in
\begin{equation*}
\mathcal{L}\xi = -\frac{1}{r^{2}}\frac{\partial V}{\partial\varphi}=\frac{1}{r^{2}}\left[4\sin(4\varphi) - \frac{\sin(0.5\varphi)}{2 \cos^2(0.5\varphi)}\right].
\end{equation*}
In order to evaluate the first term in Eq.~\eqref{eq:projected_drift_diff}, we note that the Jacobian determinant of $\xi$ is $J(x,y)=\frac{1}{r^{2}}$, and that the stationary distribution factors, canceling the $\varphi$-dependent term. We are left with
\begin{align*}
b^{\xi}(\varphi) & = \frac{1}{\vartheta(\varphi)}\int_{0}^{\infty}\frac{1}{r^{2}}\left[4\sin(4\varphi) - \frac{\sin(0.5\varphi)}{2 \cos^2(0.5\varphi)}\right] \,r\mu(r,\varphi)\,\mathrm{d}r\\
 & = \frac{1}{C_2}\left[4\sin(4\varphi) - \frac{\sin(0.5\varphi)}{2 \cos^2(0.5\varphi)}\right] \int_{0}^{\infty}\frac{1}{r}\exp(-10(r-1)^{2}-\frac{1}{r})\,\mathrm{d}r\\
 & = \frac{C_1}{C_2}\left[4\sin(4\varphi) - \frac{\sin(0.5\varphi)}{2 \cos^2(0.5\varphi)}\right],
\end{align*}
with the definitions $C_1 := \int_0^\infty \frac{1}{r}\exp(-10(r-1)^2 - \frac{1}{r})\,\mathrm{d}r$ and $C_2 := \int_0^\infty r\exp(-10(r-1)^2 - \frac{1}{r})\,\mathrm{d}r$. For the diffusion, we obtain

\begin{align*}
\nabla\xi^{T}a\nabla\xi & = 2\left[\left(\frac{\partial\xi}{\partial x}\right)^{2}+\left(\frac{\partial\xi}{\partial y}\right)^{2}\right]\\
 & = 2J(x,y)=\frac{2}{r^{2}}.
\end{align*}
Inserting this into the second term in \eqref{eq:projected_drift_diff}, and using the same arguments as before, we end up with
\begin{equation*}
a^\xi(\varphi) = \frac{2C_1}{C_2}.
\end{equation*}

\section{Parameters of Prototypical Molecular Example}
\label{sec:appendix_molecule}
The system we considered in section \ref{ssec:toy_molecule} is defined as follows: denote the three-dimensional position vectors of atoms one through five by $\mathbf{r}_i$, $1 \leq i \leq 5$. We denote the pairwise distances between atoms $i,\, j$ by $d_{i,j}$, the bond angle formed by atoms $i, j, k$ by $\theta_{i,j,k}$, and the dihedral angle formed by atoms $i, j, k, l$ by $\phi_{i,j,k,l}$. The potential energy is then defined as a superposition of harmonic bond terms, harmonic bond angle terms, and dihedral angle terms:

\begin{equation}
\label{eq:potential_toy_molecule}	
\begin{aligned}
V(\mathbf{r}_1, \ldots, \mathbf{r}_5) &= \sum_{i = 1}^4 \frac{1}{2} \,k_b (d_{i,i+1} - d^0)^2 + \\
& \quad + \sum_{i = 1}^3 \frac{1}{2}\, k_\theta (\theta_{i, i+1, i+2} - \theta^0)^2 \\
& \quad + \sum_{i=1}^2 k_\phi  (1 - \cos(n_\phi^i \,  \phi_{i,i+1,i+2,i+3})).
\end{aligned}
\end{equation}
The constants appearing above are set to
\begin{align*}
k_b &= 1, & d^0 &= 2, & k_\theta &= 1, & \theta^0 &= \frac{\pi}{2}, & k_\phi &= 0.02, & n_\phi^1 &= 3, & n_\phi^2 &= 2.
\end{align*}
The overdamped Langevin dynamics \eqref{eq:overdamped_Langevin} with potential \eqref{eq:potential_toy_molecule} are simulated at inverse temperature $\beta = 15$, for $10^7$ discrete steps at integration time step $\Delta_t = 10^{-3}$. For the purpose of our analysis, we retain only every fifth step of this trajectory.

\bibliographystyle{unsrt}
\bibliography{library}

\end{document}